\documentclass[12pt,leqno]{article}
\usepackage{amssymb,amsthm,amsmath,bbm}
\usepackage{equation}
\def\bm#1{\mathbbm{#1}}
\def\fn#1{\mathop{{\rm #1}\vphantom{\dim}}}
\def\ul#1{\mathop{\underline{#1}}}

\makeatletter
\renewcommand{\section}{\@startsection{section}{1}{0mm}{12mm}{5mm}{\raggedright\bf\Large}}

\def\@citex[#1]#2{\if@filesw\immediate\write\@auxout{\string\citation{#2}}\fi
  \def\@citea{}\@cite{\@for\@citeb:=#2\do
    {\@citea\def\@citea{\@citesep}\@ifundefined
       {b@\@citeb}{{\bf ?}\@warning
       {Citation `\@citeb' on page \thepage \space undefined}}%
{\csname b@\@citeb\endcsname}}}{#1}}
\def\@citesep{; }
\makeatother

\newtheorem{theorem}{\indent Theorem}[section]

\newtheorem{defn}[theorem]{\indent Definition}
\newtheorem{lemma}[theorem]{\indent Lemma}
\newtheorem{coro}[theorem]{\indent Corollary}
\newtheorem{proposition}[theorem]{\indent Proposition}
\newtheorem{remark}{\indent Remark}[section]

\newenvironment{idef}[1]{\par\smallskip\textit{#1}}{}

\title{Rationality of Three-Dimensional Quotients by Monomial Actions}
\author{Ming-chang Kang$^{(1)}$ and Yuri \ G.\ Prokhorov$^{(2)}$}
\date{}

\begin{document}

\maketitle
\begin{description}
\item[] $^{(1)}$Department of Mathematics and Taida Institute of
Mathematical Sciences, National Taiwan University, Taipei, Taiwan;
E-mail: kang@math.ntu.edu.tw \item[] $^{(2)}$Department of
Algebra, Faculty of Mathematics and Mechanics, Moscow State
University, Moscow, Russia; E-mail: prokhorov@gmail.com
\end{description}

\footnote{\vspace*{-5mm}\begin{description} \item[] The
first-named author was partially supported by National Center for
Theoretic Sciences (Taipei office). Parts of the work in this
paper were done while he visited Department of Mathematics, Yunnan
Normal University, Kunming. \\[-6mm]
\item[] The second-named author was partially supported by the
Russian Foundation for Basic Research (grants No~
06-01-72017-MNTI\-a, 08-01-00395-a) and Leading Scientific Schools
(grants No~ NSh-1983.2008.1, NSh-1987.2008.1).
\end{description} }

\begin{abstract}
\noindent Abstract. Let $G$ be a finite 2-group and $K$ be a field
satisfying that (i) $\fn{char}K\ne 2$, and (ii) $\sqrt{a}\in K$
for any $a\in K$. If $G$ acts on the rational function field
$K(x,y,z)$ by monomial $K$-automorphisms, then the fixed field
$K(x,y,z)^G$ is rational (= purely transcendental) over $K$.
Applications of this theorem will be given.
\end{abstract}

\bigskip


\newpage
\section{Introduction}

Let $K$ be any field and $K(x_1,\ldots,x_n)$ be the rational function field of $n$ variables over $K$.
A $K$-automorphism $\sigma$ of $K(x_1,\ldots,x_n)$ is said to be a monomial automorphism if
\[
\sigma(x_j)=a_j(\sigma) \prod_{i=1}^n x_i^{m_{ij}}
\]
where $(m_{ij})_{1\le i,j\le n} \in GL_n(\bm{Z})$ and
$a_j(\sigma)\in K\backslash \{0\}$. If $a_j(\sigma)=1$ for all
$1\le j\le n$, then $\sigma$ is called a purely monomial
automorphism. A group action on $K(x_1,\ldots,x_n)$ by monomial
$K$-automorphisms is also called a multiplicative group action.
Multiplicative actions are crucial in solving rationality problems
for linear group actions, e.g. Noether problem (see, for examples,
\cite{Sw,CHKP}).

However, the rationality problem of a multiplicative action on $K(x_1,\ldots,x_n)$ is rather intricate.
We recall some previously known results first.

\begin{theorem}[Hajja \cite{Ha1,Ha2}] \label{t1.1}
Let $K$ be any field and $K(x,y)$ be the rational function field of two variables over $K$.
If $G$ is a finite group acting on the rational function field $K(x,y)$ by monomial $K$-automorphisms,
then the fixed field $K(x,y)^G$ is rational $(=$ purely transcendental$)$ over $K$.
\end{theorem}

\begin{theorem}[Hajja, Kang, Hoshi and Rikuna \cite{HK1,HK2,HR}] \label{t1.2}
Let $K$ be any field and $K(x,y,z)$ be the rational function field of three variables over $K$.
If $G$ is a finite group acting on $K(x,y,z)$ by purely monomial $K$-automorphisms,
then the fixed field $K(x,y,z)^G$ is rational over $K$.
\end{theorem}

In the above Theorem \ref{t1.2},
it is impossible to replace the assumption ``purely monomial $K$-automorphisms" by the weaker assumption
``monomial $K$-automorphisms" as the fixed field $K(x,y,z)^{\langle\sigma\rangle}$ is not rational over $K$
where the monomial $K$-automorphism $\sigma$ is defined by $\sigma:x \mapsto y \mapsto z\mapsto -1/(xyz)$
(see the last paragraph of \cite{Ha1}).
However, if we assume that $K$ is an algebraically closed field,
then we can get an affirmative result for monomial group actions.
In fact, the main result of this article is the following theorem.

\begin{theorem} \label{t1.3}
Let $G$ be a finite $2$-group and $K$ be a field satisfying that
\begin{enumerate}
\item[{\rm (i)}] $\fn{char}K\ne 2$, and \item[{\rm (ii)}]
$\sqrt{a}\in K$ for any $a\in K$.
\end{enumerate}

If $G$ acts on the rational function field $K(x,y,z)$ by monomial
$K$-automorphisms, then the fixed field $K(x,y,z)^G$ is rational
over $K$.
\end{theorem}

In the above theorem, it is assumed that $\fn{char}K\ne 2$. We may
as well show that $K(x,y,z)^G$ is rational over $K$ when
$\fn{char}K=2$ and $K$ satisfies the assumption that $\sqrt{a}\in
K$ for any $a\in K$, by using the techniques in \cite{HK1,HK2} and
the method developed in Section 3 of this paper. Because we intend
to highlight the main techniques in the proof of Theorem
\ref{t1.3} and because we try to shorten the length of this
article, we omit the proof of the case when $\fn{char}K=2$.

We will give another remark on Theorem \ref{t1.3}. It is
impossible to generalize Theorem \ref{t1.3} to the case of
rational function fields of four variables. In fact, if $G \simeq
C_2^3$, there is a monomial action of $G$ on
$\bm{C}(x_1,x_2,x_3,x_4)$ so that the fixed field
$\bm{C}(x_1,x_2,x_3,x_4)^G$ is not retract rational over $\bm{C}$;
in particular, it is not rational over $\bm{C}$ \cite[Example
5.11]{CHKK}.

An application of Theorem \ref{t1.3} is the following theorem.
\begin{theorem} \label{t1.4}
Let $G$ be a finite $2$-group and $K$ be a field satisfying that
\begin{enumerate}
\item[{\rm (i)}] $\fn{char}K\ne 2$, and \item[{\rm (ii)}]
$\sqrt{a}\in K$ for any $a\in K$.
\end{enumerate}

Suppose that $\rho: G\to GL(V)$ is a linear representation of $G$
over $K$ such that either $\dim_K V\le 5$ or $\rho$ is the direct
sum of three $2$-dimensional irreducible representations. Then
both the fixed field $K(V)^G$ and the quotient $P(V)/G$ are
rational over $K$.
\end{theorem}

We remark that, if we assume that the field $K$ is algebraically
closed, the conclusion of the above Theorem \ref{t1.3} is still
valid for any finite group $G$, i.e. no assumption about $G$ being
a $2$-groups is necessary. We refer the reader to \cite[Section
5]{Pr} for some key ideas of the proof. We choose to publish parts
of the general result because of two reasons. First of all, the
proof of the general case is rather long and complicated; we had
better publish this result in two separate papers. Second, the
result for 2-groups will be used in a forthcoming paper on
Noether's problem and Bogomolov multipliers (i.e. the unramified
Brauer groups) for groups of order 64 \cite{CHKK}.

We will point out that Theorem \ref{t1.4} provides a quick proof
for the special case when $K$ is algebraically closed of the
following theorem.

\begin{theorem}[Chu, Hu, Kang and Prokhorov \cite{CHKP}] \label{t1.5}
Let $G$ be a group of order $32$ with exponent $e$, and $K$ be a
field satisfying

{\rm (i)} $\fn{char}K=2$, or

{\rm (ii)} $\fn{char}K\ne 2$ and $K$ contains a primitive $e$-th
root of unity,

Then $K(x_g:g\in G)^G$ is rational over $K$.
\end{theorem}

Note that Theorem \ref{t1.5} was proved using the classification
of groups of order $32$, i.e. the structures and representations
of these groups provided by the data base of GAP \cite{CHKP}. If
we choose to avoid using the classification of groups of order
$32$, we still have a proof, thanks to Theorem \ref{t1.4}, but we
need the stronger assumptions that $\fn{char}K\ne 2$ and
$\sqrt{a}\in K$ for any $a\in K$ (e.g. $K$ is algebraically
closed).

We will organize this paper as follows. In Section 2 we recall
several preliminary results which will be used subsequently.
Section 3 contains the complete proof of Theorem \ref{t1.3}; the
proof uses only elementary methods and will be more accessible to
most readers. Another proof of Theorem \ref{t1.3} for the case
when the ground field $K$ is algebraically closed and $\fn{char}K
= 0$ will be given in Section 4. This new proof is shorter than
that in Section 3, more conceptual, and no classification of
finite subgroups of $GL_3(\bm{Z})$ \cite{Ta} is required; but the
price is that some machinery in algebraic geometry is used. The
proof of Theorem \ref{t1.4} will be given in Section 5. As an
application of it, we will give a quick proof of Theorem
\ref{t1.5} provided that $K$ is algebraically closed.

Standing Notations.
Throughout this paper, $K$ is a field,
$K(x_1,\ldots,x_n)$ or $K(x,y,z)$ is the rational function field over $K$.
We will denote by $\zeta_n$ a primitive $n$-th root of unity.
Whenever we write $\zeta_n\in K$,
it is understood that either $\fn{char}K=0$ or $\fn{char}K>0$ with $\fn{char}K \nmid n$.

If $G$ is a finite group, the exponent of $G$ is $\fn{lcm} \{\fn{ord}(g): g\in G\}$
where $\fn{ord}(g)$ is the order of an element $g\in G$.
The cyclic group of order $n$ and the dihedral group of order $2n$ will be denoted by $C_n$ and $D_n$ respectively.

If $G$ is a finite group acting on $K(x_1,\ldots,x_n)$ by
$K$-automorphism, the actions of $G$ are called monomial actions
if for any $\sigma\in G$, any $1\le j\le n$, $\sigma\cdot
x_j=a_j(\sigma)\cdot \prod_{1\le i\le n}x_i^{m_{ij}}$ where
$m_{ij}\in\bm{Z}$ and $a_j(\sigma)\in K\backslash \{0\}$. The
actions are called purely monomial actions if they are monomial
actions satisfying $a_j(\sigma)=1$ for any $\sigma\in G$, any
$1\le j\le n$. In the rational function field $K(x_1,\ldots,x_n)$,
an element of the form $x_1^{\lambda_1} x_2^{\lambda_2} \cdots
x_n^{\lambda_n}$ where each $\lambda_i\in\bm{Z}$ will be called a
power product in $x_1,\ldots,x_n$; we refrain from calling it a
monomial to avoid possible confusion with the monomial action and
with $a\cdot x_1^{\lambda_1}\cdots x_n^{\lambda_n}$ where $a\in
K$.

\section{Preliminaries}

In this section we recall several results which will be used in the proof of Theorem \ref{t1.3}.

\begin{theorem}[{\cite[Theorem 1]{HK3}}] \label{t2.1}
Let $G$ be a finite group acting on $L(x_1,\ldots,x_n)$,
the rational function field of $n$ variables over a field $L$.
Suppose that

{\rm (i)} for any $\sigma\in G$, $\sigma(L)\subset L$,

{\rm (ii)} the restriction of the action of $G$ to $L$ is
faithful, and

{\rm (iii)} for any $\sigma \in G$,
\[
\begin{pmatrix} \sigma(x_1) \\ \sigma(x_2) \\ \vdots \\ \sigma(x_n) \end{pmatrix}
=A(\sigma)\cdot \begin{pmatrix} x_1 \\ x_2 \\ \vdots \\ x_n \end{pmatrix} +B(\sigma)
\]
where $A(\sigma)\in GL_n(L)$ and $B(\sigma)$ is an $n\times 1$ matrix over $L$.

Then there exist elements $z_1,\ldots,z_n\in L(x_1,\ldots,x_n)$
which are algebraically independent over $L$ and
$L(x_1,\ldots,x_n)=L(z_1,\ldots,z_n)$ so that $\sigma(z_i)=z_i$
for any $\sigma\in G$, any $1\le i\le n$.
\end{theorem}

\begin{theorem}[{\cite[Theorem 3.1]{AHK}}] \label{t2.2}
Let $L$ be any field, $L(x)$ the rational function field of one
variable over $L$, and $G$ a finite group acting on $L(x)$.
Suppose that, for any $\sigma\in G$, $\sigma(L)\subset L$ and
$\sigma(x)=a_\sigma\cdot x+b_\sigma$ where $a_\sigma,b_\sigma\in
L$ and $a_\sigma \ne 0$. Then $L(x)^G=L^G(f)$ for some polynomial
$f\in L[x]$. In fact, if $m=\min\{\deg g(x): g(x)\in
L[x]^G\setminus L^G \}$, any polynomial $f\in L[x]^G$ with $\deg
f=m$ satisfies the property $L(x)^G=L^G(f)$.
\end{theorem}

\begin{theorem}[{\cite[Theorem 2.3]{CHK}}] \label{t2.3}
Let $K$ be any field, $a,b\in K\backslash \{0\}$ and $\sigma:K(x,y)\to K(x,y)$ be a $K$-automorphism of the
rational function field $K(x,y)$ defined by $\sigma(x)=a/x$, $\sigma(y)=b/y$.
Then $K(x,y)^{\langle\sigma\rangle}=K(u,v)$ where
\[
u=\frac{x-\frac{a}{x}}{xy-\frac{ab}{xy}},~~~ v=\frac{y-\frac{b}{y}}{xy-\frac{ab}{xy}}.
\]

Moreover,
\begin{eqnarray*}
x+(a/x)=(-bu^2+av^2+1)/v, \\y+(b/y)=(bu^2-av^2+1)/u,
\\xy+(ab/(xy))=(-bu^2-av^2+1)/(uv),
\end{eqnarray*}

and

\begin{equation}
\frac{x-\frac{a}{x}}{\frac{bx}{y}-\frac{ay}{x}} =\frac{u}{bu^2-av^2}. \label{eq2.1}
\end{equation}
\end{theorem}

\begin{proof}
The proof of Formula \eqref{eq2.1} was given in \cite[p.156]{CHK}.
\end{proof}

\begin{theorem}[{\cite[Theorem 2.4]{Ka1}}] \label{t2.4}
Let $K$ be any field. Define a $K$-automorphism $\sigma$ on the
rational function field $K(x,y)$ by $\sigma(x)=a/x$, $\sigma(y)=\{
b_1[x+(a/x)]+b_2 \}/y$ where $a,b_1,b_2\in K$ with $ab_1\ne 0$.
 Then
$K(x,y)^{\langle\sigma\rangle}=K(u,v)$ where
\[
u=\frac{x-\frac{a}{x}}{xy-\frac{ab}{xy}},~~~ v=\frac{y-\frac{b}{y}}{xy-\frac{ab}{xy}}.
\]
with $b=b_1[x+(a/x)]+b_2$.
\end{theorem}

\begin{theorem}[{\cite[Theorem 6.7; HK2, Theorem 2.4]{HKO}}] \label{t2.5}
Let $F$ be a field with $\fn{char}F\ne 2$, $E=F(\alpha)$ be a field extension of $F$ defined by
$\alpha^2=a\in F\backslash \{0\}$.
Let $E(x,y)$ be the rational function field of two variables over $E$ and $\sigma$ be an $F$-automorphism on $E(x,y)$
defined by
\[
\sigma(\alpha)=-\alpha,~~\sigma(x)=x,~~\sigma(y)=b(x^2-c)/y
\]
where $b,c\in F$ with $b\ne 0$. Then
$E(x,y)^{\langle\sigma\rangle}$ is rational over $F$ if and only
if the Hilbert symbol $(a,b)_2$ is trivial in the field
$F(\sqrt{ac})$.
\end{theorem}

\begin{theorem}[{\cite[Theorem 1.4]{Ka4}}] \label{t2.6}
Let $K$ be a field and $G$ be a finite group. Assume that

{\rm (i)} $G$ contains an abelian normal subgroup $H$ so that
$G/H$ is cyclic of order $n$;

{\rm (ii)} $\bm{Z}[\zeta_n]$ is a unique factorization domain; and

{\rm (iii)} $\zeta_e\in K$ where $e$ is the exponent of $G$.

If $G\to GL(V)$ is any finite-dimensional linear representation of
$G$ over $K$, then $K(V)^G$ is rational over $K$.
\end{theorem}

\begin{defn} \label{d2.7}
Let $G$ be a finite group acting on $K(x_1,\ldots,x_n)$ by
monomial $K$-automorphisms. Let $z_1,z_2,\ldots,z_n\in
K(x_1,\ldots,x_n)$ satisfy

{\rm (i)} $\fn{trdeg}_K K(z_1,\ldots,z_n)=n$,

{\rm (ii)} for any $\sigma\in G$, any $1\le j\le n$, $\sigma\cdot
z_j=c_j(\sigma)\cdot \prod_{1\le i\le n}z_i^{m_{ij}}$ where
$m_{ij}\in\bm{Z}$ and $c_j(\sigma)\in K\backslash \{0\}$.

We will define a group homomorphism $\rho_{\ul{z}}:G\to
GL_n(\bm{Z})$ by $\rho_{\ul{z}}(\sigma)=(m_{ij})_{1\le i,j\le n}$
if $\sigma\cdot z_j$ is given by {\rm (ii)}. Note that $\ul{z}$
denotes the ordered transcendental basis $(z_1,\ldots,z_n)${\rm ;}
if $\ul{z}$ is understood from the context, we will simply write
$\rho$ for $\rho_{\ul{z}}$.
\end{defn}

\begin{lemma} \label{l2.8}
Let $K$ be any field and $G$ be a finite group acting on $K(x_1,\ldots,x_n)$ by monomial $K$-automorphisms.
Then there is a normal subgroup $H$ of $G$ so that

{\rm (i)} $K(x_1,\ldots,x_n)^H=K(z_1,\ldots,z_n)$ for some power
products $z_1,\ldots,z_n$ in $x_1,\ldots,x_n$,

{\rm (ii)} $G/H$ acts on $K(z_1,\ldots,z_n)$ by monomial
$K$-automorphisms, and

{\rm (iii)} $\rho_{\ul{z}}: G/H \to GL_n(\bm{Z})$ is injective.
\end{lemma}

\begin{proof}
Induction on the order of $G$.

Without loss of generality, we may assume that $G$ acts faithfully on $K(x_1,\ldots,x_n)$.
Define $H_0=\fn{Ker} \{\rho_{\ul{x}}: G\to GL_n(\bm{Z})\}$.
For any $\tau\in H_0$, $\tau(x_i)=a_i(\tau)\cdot x_i$ for some $a_i(\tau)\in K\backslash \{0\}$.
In particular, $H_0$ is an abelian subgroup of $G$.
Choose $\tau_1,\ldots,\tau_m\in H_0$ so that $H_0$ is generated by $\tau_1,\ldots,\tau_m$.

Define $\langle x_1,\ldots,x_n\rangle$ to be the multiplicative
subgroup of $K(x_1,\ldots,x_n)\backslash\{0\}$ generated by
$x_1,\ldots,x_n$, i.e. $\langle
x_1,\ldots,x_n\rangle=\{M=x_1^{\lambda_1} x_2^{\lambda_2} \cdots
x_n^{\lambda_n}: \lambda_i\in\bm{Z}$ for $1\le i\le n\}$. Define a
group homomorphism
\begin{eqnarray*}
\Phi: \langle x_1,\ldots,x_n\rangle &\to& K^{\times}\times K^{\times} \times \cdots \times K^{\times} \\
M &\mapsto& (\tau_1(M)/M, \tau_2(M)/M,\ldots,\tau_m(M)/M)
\end{eqnarray*}
where $K^{\times}=K\backslash \{0\}$ and $M=x_1^{\lambda_1} \cdots
x_n^{\lambda_n}$ with $\lambda_i\in\bm{Z}$. Note that the image of
$\Phi$ is a finite group. It is not difficult to find that
$\fn{Ker}(\Phi)=\langle M_1,\ldots,M_n\rangle$ where
$M_1,\ldots,M_n$ are $n$ power products in $x_1,\ldots,x_n$.
Moreover, $K(x_1,\ldots,x_n)^{H_0}=K(M_1,\ldots,M_n)$.

We will show that $G$ (equivalently, $G/H_0$) acts on
$K(M_1,\ldots,M_n)$ by monomial $K$-automorphisms. For any
$\sigma\in G$, any $1\le i\le n$, $\sigma(M_i)=c_i(\sigma)\cdot
N_i$ where $N_i=x_1^{v_1} \cdots x_n^{v_n}$ with $v_i\in \bm{Z}$
and $c_i(\sigma)\in K\backslash \{0\}$. To show that each $N_i$ is
a power product in $M_1,M_2,\ldots,M_n$, it suffices to show that
$N_i\in K(x_1,\ldots,x_n)^{H_0}=K(M_1,\ldots,M_n)$. From the
identity $\sigma(M_i)=c_i(\sigma)\cdot N_i$, we find that, for any
$\tau\in H_0$, $c_i(\sigma)\cdot
\tau(N_i)=\tau\sigma(M_i)=\sigma\cdot(\sigma^{-1}\tau\sigma)(M_i)=\sigma(M_i)$
because $\sigma^{-1}\tau\sigma\in H_0$. Thus $\tau(N_i)=N_i$ for
any $\tau\in H_0$.

Now consider $\rho_{\ul{M}}:G/H_0 \to GL_n(\bm{Z})$. By induction
hypothesis, we can find $z_1$, $\ldots$, $z_n\in
K(M_1,\ldots,M_n)$ so that (i) $z_1,\ldots,z_n$ are power products
in $M_1,\ldots,M_n$, (ii)
$K(z_1,\ldots,z_n)=K(M_1,\ldots,M_n)^{H/H_0}$ for some normal
subgroup $H/H_0$ in $G/H_0$ (where $H \supset H_0$), (iii) $G/H$
acts monomially on $K(z_1,\ldots,z_n)$, (iv) $\rho_{\ul{z}}: G/H
\to GL_n(\bm{Z})$ is injective.
\end{proof}

\setcounter{equation}{0}
\section{Proof of Theorem \ref{t1.3}}

We will prove Theorem \ref{t1.3} in this section. Because of Lemma
\ref{l2.8}, it suffices to consider the situation when $\rho: G\to
GL_3(\bm{Z})$ is injective. This condition will remain in force
throughout this section.

Unless otherwise specified, we will assume, in this section, that
$G$ is a finite 2-group acting on $K(x_1,x_2,x_3)$ by monomial
$K$-automorphisms, and $K$ is a field satisfying the conditions
(i) $\fn{char} K\ne 2$, and (ii) $\sqrt{a}\in K$ for any $a\in K$.
For example, in the following Theorem \ref{t3.4} it is understood
that $K$ satisfies the above two conditions although there is no
mention about the assumptions on $K$ in the statement of Theorem
\ref{t3.4}. On the other hand, in Theorem \ref{t3.3}, $K$
satisfies the weaker conditions that $\fn{char}K\ne 2$ and
$\sqrt{-1}\in K$, which will be stated explicitly.

Since $\rho:G\to GL_3(\bm{Z})$ is injective, $G$ is isomorphic to
a finite subgroup of $GL_3(\bm{Z})$ as an abstract group. It is
known that, up to conjugation in $GL_3(\bm{Z})$, there are
precisely 73 non-isomorphic finite subgroups in $GL_3(\bm{Z})$
\cite[p.807]{Ta,HK2}. We will denote by $W_i(j)$ the group $W_i$
which appears on page $j$ of Tahara's paper \cite{Ta}. The
2-groups in $GL_3(\bm{Z})$ are the following 36 groups $G$,
\leftmargini=18mm
\begin{enumerate}
\item[(I)] $G\simeq C_2: W_i(173)$ where $1\le i\le 5$;
\item[(II)] $G\simeq C_4: W_i(174)$ where $1\le i\le 4$; and \\
$G\simeq C_2\times C_2:W_i(174)$ where $5\le i\le 15$;
\item[(III)] $G\simeq C_4\times C_2: W_i(187)$ where $1\le i\le 2$; \\
$G\simeq C_2\times C_2\times C_2:W_i(187)$ where $3\le i\le 6$; and \\
$G\simeq D_4: W_i(187)$ where $7\le i\le 14$;
\item[(IV)] $G\simeq D_4\times C_2: W_1(194)$ and $W_2(195)$.
\end{enumerate}

Convention.
Suppose that $\sigma,\tau\in G$.
We will adopt the convention
\begin{eqnarray*}
\sigma &:& (x_1,x_2,x_3) \mapsto (\epsilon_1x_1,1/x_2,1/x_3), \\
\tau &:& (x_1,x_2,x_3) \mapsto (1/x_1,\epsilon_2/x_2,\epsilon_3/x_3),
\end{eqnarray*}
to indicate the fact that $\sigma$ and $\tau$ are $K$-automorphisms on $K(x_1,x_2,x_3)$ defined by
$\sigma(x_1)=\epsilon_1x_1$, $\sigma(x_2)=1/x_2$, $\sigma(x_3)=1/x_3$,
$\tau(x_1)=1/x_1$, $\tau(x_2)=\epsilon_2/x_2$, $\tau(x_3)=\epsilon_3/x_3$.

\begin{theorem} \label{t3.1}
Let $G=W_i(173)$ for $1\le i\le 4$ and $K$ be a field with $\fn{char}K\ne 2$.
Then $K(x_1,x_2,x_3)^G$ is rational over $K$.
\end{theorem}

\begin{proof}
First we will ``normalize" the coefficients of $\sigma(x_i)$ where $G=\langle\sigma\rangle$.

When $G=W_4(173)$, the action of $G=\langle\sigma\rangle$ is given by
\[
\sigma: (x_1,x_2,x_3) \mapsto (ax_1,b/x_3,c/x_2)
\]
for some $a,b,c\in K \backslash \{0\}$. Since $G\simeq C_2$, it
follows that $a=\pm 1$. Define $y_2=x_2$, $y_3=b/x_3$. Then
$K(x_2,x_3)=K(y_2,y_3)$ and $\sigma:(y_2,y_3) \mapsto
(y_3,by_2/c)$. Since $G\simeq C_2$, we find that $b/c=1$. Apply
Theorem \ref{t2.1}. We find that
$K(x_1,x_2,x_3)^G=K(y_2,y_3,x_1)^G=K(y_2,y_3)^G(u)$ for some $u$
with $\sigma(u)=u$. Note that
$K(y_2,y_3)^G=K(y_2+y_3,(y_2-y_3)^2)$ is rational over $K$.

When $G=W_1(173)$, the action of $G$ is given by
\[
(x_1,x_2,x_3)\mapsto(ax_1, b/x_2, c/x_3) \] for some $a,b,c\in
K\backslash\{0\}$. Apply Theorem \ref{t2.1} and reduce the
question to $K(x_2,x_3)^G$. Apply Theorem \ref{t1.1}. Done.

When $G=W_2(173)$ and $W_3(173)$, the actions of $G$ are given by
\[
(x_1,x_2,x_3)\mapsto (a/x_1, bx_2, cx_3) \] and
\[
(x_1,x_2,x_3)\mapsto (a/x_1, bx_3, cx_2) \] respectively where
$a,b,c\in K\backslash\{0\}$.

Apply Theorem \ref{t2.1} to both cases. The questions are reduced
to $K(x_1)^G$. Note that $K(x_1)^G$ is rational by L\"{u}roth's
Theorem.

\end{proof}

\begin{theorem} \label{t3.2}
Let $K$ be a field with $\fn{char}K\ne 2$ and $G=W_5(173)$,
i.e. the action of $G=\langle\sigma\rangle$ is given by
\[
\sigma: (x_1,x_2,x_3) \mapsto (a_1/x_1,a_2/x_2,a_3/x_3)
\]
for some $a_1,a_2,a_3\in K\backslash\{0\}$.
Then $K(x_1,x_2,x_3)^G$ is rational over $K$ if and only if $[K(\sqrt{a_1},\sqrt{a_2},\sqrt{a_3}):K]\le 4$.
\end{theorem}

\begin{proof}
If $[K(\sqrt{a_1},\sqrt{a_2},\sqrt{a_3}):K]=8$, it is Saltman who
shows that $K(x_1,x_2,x_3)^G$ is not retract rational over $K$
\cite{Sa} (note that $K$ is an infinite field in this situation).
In particular, $K(x_1,x_2,x_3)^G$ is not rational over $K$. See
\cite{Ka2} for a generalization of Saltman's Theorem.

If $[K(\sqrt{a_1},\sqrt{a_2},\sqrt{a_3}):K] \le 4$, we may assume
that $a_3\in a_1^{\epsilon_1}a_2^{\epsilon_2}K^2$ for some
$\epsilon_1,\epsilon_2\in \{0,1\}$. First consider the case
$a_3=a_1a_2b^2$ for some $b\in K\backslash\{0\}$. Define
$y=x_3/(bx_1x_2)$. Then $K(x_1,x_2,x_3)=K(x_1,x_2,y)$ and
$\sigma(y)=1/y$. Define $z=(1-y)/(1+y)$. Then $\sigma(z)=-z$. It
follows that $K(x_1,x_2,x_3)^G=K(x_1,x_2,z)^G=K(x_1,x_2)^G(u)$ for
some $u$ with $\sigma(u)=u$ by Theorem \ref{t2.1}. Note that
$K(x_1,x_2)^G$ is rational by Theorem \ref{t1.1}. The other cases
$a_3\in a_1K^2$, $a_3\in a_2K^2$ and $a_3\in K^2$ can be proved
similarly.
\end{proof}

\begin{theorem} \label{t3.3}
Let $G=W_i(174)$ for $1\le i\le 4$ and $K$ satisfy the conditions that $\fn{char}K\ne 2$ and $\sqrt{-1}\in K$.
Then $K(x_1,x_2,x_3)^G$ is rational over $K$.
\end{theorem}
\begin{proof}
\begin{idef}{Case 1.} $G=W_1(174)$.

The action of $G=\langle\sigma\rangle$ is given by
\[
\sigma: (x_1,x_2,x_3) \mapsto \left( (\sqrt{-1})^ix_1,x_3,a/x_2\right)
\]
for some integer $i$ and some $a\in K\backslash\{0\}$.
Apply Theorem \ref{t2.1}.
We find that $K(x_1$, $x_2,x_3)^G=K(x_2,x_3)^G(u)$ for some $u$ with $\sigma(u)=u$.
By Theorem \ref{t1.1}, $K(x_2,x_3)^G$ is rational over $K$.
\end{idef}

\bigskip
\begin{idef}{Case 2.} $G=W_2(174)$.

The action of $G=\langle\sigma\rangle$ is given by
\[
\sigma: (x_1,x_2,x_3) \mapsto (a/x_1,x_3,b/x_2)
\]
for some $a,b\in K\backslash\{0\}$.
Define $u$ and $v$ by
\[
u=\frac{x_2-\frac{b}{x_2}}{x_2x_3-\frac{b^2}{x_2x_3}}, ~~~
v=\frac{x_3-\frac{b}{x_3}}{x_2x_3-\frac{b^2}{x_2x_3}}.
\]

By Theorem \ref{t2.3},
we get $K(x_1,x_2,x_3)^{\langle\sigma^2\rangle}=K(x_1,u,v)$.
Note that
\[
\sigma: u\mapsto
\frac{x_3-\frac{b}{x_3}}{\frac{bx_3}{x_2}-\frac{bx_2}{x_3}}, ~ v
\mapsto
\frac{-\big(x_2-\frac{b}{x_2}\big)}{\frac{bx_3}{x_2}-\frac{bx_2}{x_3}}.
\]

Note that
\[
\frac{-\big(x_2-\frac{b}{x_2}\big)}{\frac{bx_3}{x_2}-\frac{bx_2}{x_3}} = \frac{u}{bu^2-bv^2}
\]
by Theorem \ref{t2.3} again.

Define $w=u/v$.
Then $K(x_1,u,v)=K(w,x_1,v)$ and
\[
\sigma: (w,x_1,v) \mapsto (-1/w,a/x_1,1/(cv))
\]
where $c=b[w-(1/w)]$.

Define
\[
y_1=x_1,~~~ Y_2=\big(\sqrt{-1}-w\big)\big/\big(\sqrt{-1}+w\big),
~~~ y_2=[x_1-(a/x_1)] Y_2.
\]

Then $K(w,x_1,v)=K(y_1,y_2,v)$ and $\sigma(y_1)=a/y_1$,
$\sigma(y_2)=y_2$, $\sigma(v)=1/(cv)$.

It is not difficult to verify that $c=b[w-(1/w)]=2\sqrt{-1}b \big(x_1-(a/x_1)+\sqrt{-1} Y_2\big) %
\big(
x_1-(a/x_1)-\sqrt{-1}Y_2\big)\big/[\big(x_1-(a/x_1)+Y_2\big)\big(x_1-(a/x_1)-Y_2\big)]$.

Define
\[
y_3=v(x_1-(a/x_1)+\sqrt{-1}Y_2)/(x_1-(a/x_1)+Y_2).
\]

It follows that $\sigma(y_3)=d/y_3$ where $d=1/(2\sqrt{-1}b)$.
Now we find that $K(w,x_1,v)^{\langle\sigma\rangle}$ $=K(y_1,y_2,y_3)^{\langle\sigma\rangle} %
=K(y_2)(y_1,y_3)^{\langle\sigma\rangle}$ is rational over $K(y_2)$ by Theorem \ref{t2.3}.
\end{idef}

\bigskip
\begin{idef}{Case 3.} $G=W_3(174)$.

The action of $G=\langle\sigma\rangle$ is given by
\[
\sigma:(x_1,x_2,x_3)\mapsto (ax_1,x_3,bx_1/x_2).
\]
for some $a,b\in K\backslash\{0\}$.
Define $y_1=x_2$, $y_2=x_3$, $y_3=bx_1/x_2$.
Then $K(x_1,x_2,x_3)=K(y_1,y_2,y_3)$ and $\sigma:(y_1,y_2,y_3)\mapsto (y_2,y_3,\epsilon y_1y_3/y_2)$.
Since $\fn{ord}(\sigma)=4$,
it follows that $\epsilon^2=1$, i.e. $\epsilon=\pm 1$.

\bigskip
\begin{idef}{Case 3.1.} $\epsilon=1$.

The action $\sigma:(y_1,y_2,y_3)\mapsto (y_2,y_3,y_1y_3/y_2)$ is a purely monomial automorphism.
Thus $K(y_1,y_2,y_3)^G$ is rational over $K$ by Theorem \ref{t1.2}.
\end{idef}

\bigskip
\begin{idef}{Case 3.2.} $\epsilon=-1$.

We find that $\sigma:(y_1,y_2,y_3) \mapsto (y_2,y_3,-y_1y_3/y_2)$.

Define
\[
\begin{array}{ll}
X_1=(y_1+y_3)/2, & X_2=[y_2-(y_1y_3/y_2)]/2, \\
X_3=(y_1-y_3)/(2\sqrt{-1}), & X_4=[y_2+(y_1y_3/y_2)]/(2\sqrt{-1}).
\end{array}
\]

Then $K(y_1,y_2,y_3)=K(X_1,X_2,X_3,X_4)$ with the relation $X_1^2+X_2^2+X_3^2+X_4^2=0$.
Note that $\sigma:X_1 \leftrightarrow X_2$, $X_3\mapsto X_4\mapsto -X_3$.
Thus $\sigma^2:X_1\mapsto X_1$, $X_2\mapsto X_2$, $X_3\mapsto -X_3$, $X_4\mapsto -X_4$.

Define $Y_3=X_3/X_4$, $Y_1=X_1Y_3$, $Y_2=X_2Y_3$, $Y_4=X_3^2$.
Then $K(X_1,X_2,X_3,X_4)^{\langle\sigma^2\rangle}=K(Y_1,Y_2,Y_3,Y_4)$ with the relation $Y_1^2+Y_2^2+Y_4(1+Y_3^2)=0$.
Hence $Y_4\in K(Y_1,Y_2,Y_3)$.
It follows that $K(y_1,y_2,y_3)^{\langle\sigma\rangle}=K(Y_1,Y_2,Y_3)^{\langle\sigma\rangle}$.
Note that
\[
\sigma: (Y_1,Y_2,Y_3) \mapsto (-Y_2/Y_3^2,-Y_1/Y_3^2,-1/Y_3).
\]

Define $z_1=Y_1/Y_3$, $z_2=Y_2/Y_3$. Then
$K(Y_1,Y_2,Y_3)=K(z_1,z_2,Y_3)$ and $\sigma:(z_1,z_2$,
$Y_3)\mapsto (z_2,z_1,-1/Y_3)$.

Apply Theorem \ref{t2.1}.
We find that $K(Y_1,Y_2,Y_3)^{\langle\sigma\rangle}=K(z_1,z_2,Y_3)^{\langle\sigma\rangle} %
=K(Y_3)^{\langle\sigma\rangle}(u_1,u_2)$ for some $u_1$, $u_2$ with $\sigma(u_1)=u_1$, $\sigma(u_2)=u_2$.
\end{idef}
\end{idef}

\bigskip
\begin{idef}{Case 4.} $G=W_4(174)$.

The action of $G=\langle\sigma\rangle$ is given by
\[
\sigma:(x_1,x_2,x_3) \mapsto (c/x_1,x_3,bx_1/x_2)
\]
for some $b,c\in K\backslash\{0\}$.

Define $y_1=x_2$, $y_2=x_3$, $y_3=bx_1/x_2$. Then
$K(x_1,x_2,x_3)=K(y_1,y_2,y_3)$ and $\sigma:(y_1,y_2,y_3) \mapsto
(y_2,y_3,a/(y_1y_2y_3))$ for some $a\in K\backslash \{0\}$.

Note that $\sigma^2: (y_1,y_2,y_3)\mapsto (y_3,a/(y_1y_2y_3),y_1)$.

Define $t=y_1y_3$. Then $K(y_1,y_2,y_3)=K(t,y_1,y_2)$ and
$\sigma^2:(t,y_1,y_2) \mapsto (t,t/y_1,a/(ty_2))$. Define $A=t$,
$B=a/t$, $u$ and $v$ as follows
\[
u=\frac{y_1-\frac{A}{y_1}}{y_1y_2-\frac{AB}{y_1y_2}}, ~~~
v=\frac{y_2-\frac{B}{y_2}}{y_1y_2-\frac{AB}{y_1y_2}}.
\]

By Theorem \ref{t2.3},
we find that $K(t,y_1,y_2)^{\langle\sigma^2\rangle}=K(t,u,v)$.
Note that
\begin{eqnarray*}
\sigma &:& t\mapsto a/t, \\
&& u\mapsto \frac{y_2-\frac{B}{y_2}}{\frac{Ay_2}{y_1}-\frac{By_1}{y_2}}, ~
v\mapsto \frac{-\big(y_1-\frac{A}{y_1}\big)}{\frac{Ay_2}{y_1}-\frac{By_1}{y_2}}.
\end{eqnarray*}

Note that
\[
\frac{-\big( y_1-\frac{A}{y_1}\big)}{\frac{Ay_2}{y_1}-\frac{By_1}{y_2}}=\frac{u}{Bu^2-Av^2}
\]
by Theorem \ref{t2.3}.

Define $w=u/v$.
Then $K(t,u,v)=K(t,w,v)$ and $\sigma:(t,w,v)\mapsto (a/t,-1/w,C/v)$ where $C=w/(Bw^2-A)$.

Define $Z_1=({\sqrt{-1}-w})/({\sqrt{-1}+w})$, $z_2=w/t$. Then
$\sigma(Z_1)=-Z_1$, $\sigma(z_2)=-1/(az_2)$.

Define $z_1=Z_1[z_2+1/(az_2)]$.
Then $\sigma:(z_1,z_2,v)\mapsto (z_1,-1/(az_2),C/v)$ where $C=w/(Bw^2-A)=1/[(aw/t)-(t/w)]=1/[az_2-(1/z_2)]$.

Define $z_3=v[az_2-(1/z_2)]$.
Then $K(t,u,v)=K(t,w,v)=K(z_1,z_2,z_3)$ and $\sigma:(z_1,z_2,z_3) \mapsto (z_1,-1/(az_2),D/z_3)$ where
$D=az_2-(1/z_2)=a[z_2+\sigma(z_2)]$.
By Theorem \ref{t2.4} (instead of Theorem \ref{t2.3} !),
$K(z_1,z_2,z_3)^{\langle\sigma\rangle}$ is rational over $K(z_1)$.
\end{idef}
\end{proof}

\begin{theorem} \label{t3.4}
Let $G=W_i(174)$ for $5\le i\le 15$.
Then $K(x_1,x_2,x_3)^G$ is rational over $K$.
\end{theorem}

\begin{proof}
\begin{idef}{Case 1.} $G=W_5(174)$.

The action of $G=\langle\sigma, \tau\rangle$ is given by
\begin{eqnarray*}
\sigma &:& (x_1,x_2,x_3) \mapsto (\epsilon x_1,b_2/x_2,b_3/x_3), \\
\tau &:& (x_1,x_2,x_3) \mapsto (a_1/x_1,a_2/x_2,a_3/x_3),
\end{eqnarray*}
for some $a_1,a_2,a_3,b_2,b_3\in K\backslash \{0\}$ and $\epsilon=\pm 1$.

Since $\sqrt{a_1},\sqrt{a_2},\sqrt{a_3}\in K$,
we may define $y_i=x_i/\sqrt{a_i}$ for $1\le i\le 3$.
It follows that $\sigma: (y_1,y_2,y_3)\mapsto(\epsilon y_1,c_2/y_2,c_3/y_3)$ and
$\tau:(y_1,y_2,y_3)\mapsto(1/y_1,1/y_2,1/y_3)$ for some $c_2,c_3\in K\backslash \{0\}$.
Since $\sigma\tau=\tau\sigma$,
it follows that $c_2^2=c_3^2=1$.

If $c_2c_3=1$, we define $y_4=y_2/y_3$.
It follows that $\sigma(y_4)=\tau(y_4)=1/y_4$.
Define $y_5=(1-y_4)/(1+y_4)$.
Then $\sigma(y_5)=\tau(y_5)=-y_5$.
Thus $K(y_1,y_2,y_3)^G=K(y_1,y_2,y_5)^G=K(y_1,y_2)^G(y_0)$ for some $y_0$ with $\sigma(y_0)=\tau(y_0)=y_0$
by Theorem \ref{t2.1}.
Now $K(y_1,y_2)^G$ is rational by Theorem \ref{t1.1}.

If $c_2c_3=-1$,
we may assume that $c_2=-1$ and $c_3=1$.
In this situation, define $y_4=y_3$.
The arguments in the preceding paragraph remain valid. Done.

It is can be shown that $K(x_1,x_2,x_3)^G$ is rational when
$G=W_5(174)$ and $K$ is any field with $\fn{char}K\ne 2$, i.e. the
assumption that $\sqrt{a}\in K$ for any $a\in K$ can be waived.
But we omit the proof here.
\end{idef}

\bigskip
\begin{idef}{Case 2.} $G=W_6(174)$.

The action of $G=\langle\sigma,\tau\rangle$ is given by
\begin{eqnarray*}
\sigma &:& (x_1,x_2,x_3) \mapsto (\epsilon_1 x_1,a_2/x_2,a_3/x_3), \\
\tau &:& (x_1,x_2,x_3) \mapsto (a_1/x_1,a_4/x_2,\epsilon_3x_3),
\end{eqnarray*}
for some $a_1,a_2,a_3,a_4\in K\backslash \{0\}$ and $\epsilon_1,\epsilon_3\in \{1,-1\}$.

Define $y_i=x_i/\sqrt{a_i}$ for $1\le i\le 3$.
Using the relation $\sigma\tau=\tau\sigma$,
we get
\begin{eqnarray*}
\sigma &:& (y_1,y_2,y_3) \mapsto (\epsilon_1y_1,1/y_2,1/y_3), \\
\tau &:& (y_1,y_2,y_3)\mapsto (1/y_1,\epsilon_2/y_2,\epsilon_3y_3)
\end{eqnarray*}
where $\epsilon_1,\epsilon_2,\epsilon_3\in \{1,-1\}$.

If at least one of $\epsilon_1,\epsilon_2,\epsilon_3$ is 1, say $\epsilon_3=1$,
define $y_4=(1-y_3)/(1+y_3)$.
Then we get $\sigma(y_3)=-y_3$ and $\sigma(y_3)=y_3$.
Then we may apply Theorem \ref{t2.1} and get $K(y_1,y_2,y_3)^G=K(y_1,y_2,y_4)^G=K(y_1,y_2)^G(y_0)$ with
$\sigma(y_0)=\tau(y_0)=y_0$.

It remains to consider the situation $\epsilon_1=\epsilon_2=\epsilon_3=-1$.

Define $z_1$, $u$ and $v$ by
\[
z_1=y_1(y_2-(1/y_2)),~~~
u=\frac{y_2-\frac{1}{y_2}}{y_2y_3-\frac{1}{y_2y_3}},~~~
v=\frac{y_3-\frac{1}{y_3}}{y_2y_3-\frac{1}{y_2y_3}}.
\]

Then $K(y_1,y_2,y_3)^{\langle\sigma\rangle}=K(z_1,u,v)$ by Theorem \ref{t2.3}.
Define $w=u/v$.
Use the same technique as in the proof of Case 2 of Theorem \ref{t3.3}.
We find that
\[
\tau:(z_1,u,w)\mapsto (A/z_1,B/u,-w)
\]
where $A=[y_2-(1/y_2)]^2=[(u^2+w^2-u^2w^2)^2/(u^2w^2)]-4$, $B=-w^2/(w^2-1)$.
Note that $\tau(A)=A$.

Define $z_2=u(w-1)/w$.
Then $K(z_1,u,v)=K(z_1,z_2,w)$ and $\tau(z_2)=-1/z_2$.
It is not difficult to verify that
\begin{eqnarray*}
A &=& z_2-(1/z_2)+z_2(1+w)^2-[(1-w)^2/z_2] \\
&=& 2w(z_2+(1/z_2))+(2+w^2)(z_2-(1/z_2)).
\end{eqnarray*}

Define $t=(\sqrt{-1}-z_2)/(\sqrt{-1}+z_2)$.
Then $K(z_1,u,w)=K(t,w,z_1)$ and
\begin{equation}
\tau: (t,w,z_1)\mapsto (-t,-w,A/z_1) \label{eq3.1}
\end{equation}
where $A=2\sqrt{-1}[w^2(1+t^2)-4wt+2(1+t^2)]/(1-t^2)$.

Define $x=tw$, $y=z_1(1-t)/e$ where $e\in K$ satisfies $e^2=2\sqrt{-1}$.
We will use Theorem \ref{t2.5} to show that $K(t,x,y)^{\langle\tau\rangle}$ is rational over $K(t^2)$;
in particular, $K(t,x,y)^{\langle\tau\rangle}$ is rational over $K$.

Note that $\tau(t)=-t$, $\tau(x)=x$, $\tau(y)=B/y$ where
\begin{eqnarray*}
B &=& w^2(1+t^2)-4wt+2(1+t^2) \\
&=& x^2[1+(1/t^2)]-4x+2(1+t^2) \\
&=& [1+(1/t^2)]\{[x-(2t^2/(1+t^2))]^2+(2t^2+2t^6)/(1+t^2)^2\}.
\end{eqnarray*}

In applying Theorem \ref{t2.5}, define $a=t^2$, $F=K(t^2)$, $b=(1+t^2)/t^2=(1+a)/a$ and $c=-(2a+2a^3)/(1+a)^2$.
It is known that $K(t,x,y)^{\langle\tau\rangle}=F(t)(x,y)^{\langle\tau\rangle}$ is rational over $F$
if and only if the Hilbert symbol $(a,b)_2$ is trivial in the field $F(\sqrt{ac})$.

We claim that $F(\sqrt{ac})$ is a rational function field of one
variable over $K$. Note that $F(\sqrt{ac})=F(\sqrt{1+a^2})$
because $\sqrt{-2}\in K\subset F$. Write $p=\sqrt{1+a^2}$,
$s=p-a$. Then $F(\sqrt{1+a^2})=K(a,p)$. From the relation
$(p-a)(p+a)=1$, we find that $p+a=1/(p-a)=1/s \in K(s)$. Thus
$a=(1-s^2)/(2s)\in K(s)$. Hence $F(\sqrt{1+a^2})=K(a,p)=K(s)$.

Now we find that $(a,b)_2=(a,(1+a)/a)_2=(a,a(1+a))_2=(a,-(1+a))_2=(a,1+a)_2$ because $-1\in K^2$.
To show that $(a,1+a)_2=0$ is equivalent to find a solution in $F(\sqrt{ac})=K(s)$ of the equation
\[
aX^2+(1+a)Y^2=1.
\]

Substitute $a=(1-s^2)/(2s)$ into the above equation. We get
\[
(1-s^2)X^2+(1+2s-s^2)Y^2=2s.
\]

Note that $X=\sqrt{-1}$, $Y=1$ is a solution we need. Hence the result.

It can be shown that the field $F(t)(x,y)^{\langle\tau\rangle}$ is
the function field of certain conic bundle of $\bm{P}^1$ defined
over $F$. Thus we may as well use Iskovskikh's theory of conic
bundles \cite{Is} to show that $F(t)(x,y)^{\langle\tau\rangle}$ is
rational over $F$. See also \cite{Ka3}.
\end{idef}

\bigskip
\begin{idef}{Case 3.} $G=W_7 (174)$.

The action of $G={\langle\sigma,\tau\rangle}$ is given by
\begin{eqnarray*}
\sigma &:& (x_1,x_2,x_3) \mapsto (\epsilon x_1,a_2/x_2,a_3/x_3) \\
\tau &:& (x_1,x_2,x_3) \mapsto (\epsilon_1
x_1,\epsilon_2x_2,b/x_3)
\end{eqnarray*}
where $a_2,a_3,b\in K\backslash \{0\}$ and $\epsilon,\epsilon_1,\epsilon_2\in \{1,-1\}$.

Apply Theorem \ref{t2.1}.
The question is reduced to the rationality of $K(x_2,x_3)^G$. Done.
\end{idef}

\bigskip
\begin{idef}{Case 4.} $G=W_8(174)$.

The action of $G={\langle\sigma,\tau\rangle}$ is given by
\begin{eqnarray*}
\sigma &:& (x_1,x_2,x_3) \mapsto (\epsilon x_1,a_2/x_2,a_3/x_3), \\
\tau &:& (x_1,x_2,x_3) \mapsto (b_1/x_1,b_2/x_3,b_3/x_2)
\end{eqnarray*}
where $\epsilon_1,a_2,a_3,b_1,b_2,b_3\in K\backslash \{0\}$.

Define $y_1=x_1/\sqrt{b_1}$, $y_2=x_2/\sqrt{a_2}$, $y_3=x_3/\sqrt{a_3}$.
We find that
\begin{eqnarray*}
\sigma &:& (y_1,y_2,y_3) \mapsto (\epsilon_1y_1,1/y_2,1/y_3), \\
\tau &:& (y_1,y_2,y_3)\mapsto (1/y_1,\epsilon_2/y_2,\epsilon_3/y_3)
\end{eqnarray*}
where $\epsilon_1,\epsilon_2,\epsilon_3\in \{1,-1\}$.

The situation is the same as in Case 1. Hence we omit the proof.
\end{idef}

\bigskip
\begin{idef}{Case 5.} $G=W_9(174)$. \

The action of $G={\langle\sigma,\tau\rangle}$ is given by
\begin{eqnarray*}
\sigma &:& (x_1,x_2,x_3) \mapsto (a_1 x_1,a_2/x_2,a_3/x_3) \\
\tau &:& (x_1,x_2,x_3) \mapsto (b_1x_1,b_2x_3,b_3x_2)
\end{eqnarray*}
where $a_i,b_j\in K\backslash \{0\}$.

Apply Theorem \ref{t2.1}. The question is reduced to the rationality of $K(x_2,x_3)^G$.
\end{idef}

\bigskip
\begin{idef}{Case 6.} $G=W_{10}(174)$.

The action of $G={\langle\sigma,\tau\rangle}$ is given by
\begin{eqnarray*}
\sigma &:& (x_1,x_2,x_3) \mapsto (a_1/x_1,a_2x_2,a_3x_3), \\
\tau &:& (x_1,x_2,x_3) \mapsto (b_1x_1,b_2x_3,b_3x_2)
\end{eqnarray*}
where $a_i,b_j\in K\backslash \{0\}$.

Define $y=x_1/\sqrt{a_1}$. Then $\sigma(y)=1/y$ and $\tau(y)=\pm y$.

If $\tau(y)=-y$, we will apply Theorem \ref{t2.1} and reduce the question to the rationality of $K(y)^G$.

If $\tau(y)=y$, define $z=(1-y)/(1+y)$. Then $\sigma(z)=-z$ and $\tau(z)=z$.
Apply Theorem \ref{t2.1} and reduced the question to the rationality of $K(x_2,x_3)^G$.
\end{idef}

\bigskip
\begin{idef}{Case 7.} $G=W_{11}(174)$.

After the ``normalization" of the coefficients, the action of $G={\langle\sigma,\tau\rangle}$ is given by
\begin{eqnarray*}
\sigma &:& (x_1,x_2,x_3) \mapsto (1/x_1,x_3,x_2), \\
\tau &:& (x_1,x_2,x_3) \mapsto (\epsilon x_1,b/x_2,b/x_3)
\end{eqnarray*}
where $b\in K\backslash \{0\}$ and $\epsilon=\pm 1$.

Define $y_i=x_i/\sqrt{b}$ for $i=2,3$. Then $\sigma:y_2
\leftrightarrow y_3$ and $\tau:y_2\mapsto 1/y_2$, $y_3\mapsto
1/y_3$.

If $\tau(x_1)=-x_1$, we may apply Theorem \ref{t2.1} and consider
the rationality of $K(x_1)^G$.

If $\tau(x_1)=x_1$, define $y_1=(1-x_1)/(1+x_1)$.
We get $\sigma(y_1)=-y_1$ and $\tau(y_1)=y_1$.
Apply Theorem \ref{t2.1} and consider the rationality of $K(y_2,y_3)^G$. Done.
\end{idef}

\bigskip
\begin{idef}{Case 8.} $G=W_{12}(174)$.

The action of $G={\langle\sigma,\tau\rangle}$ is given by
\begin{eqnarray*}
\sigma &:& (x_1,x_2,x_3) \mapsto (1/x_1,x_3,x_2), \\
\tau &:& (x_1,x_2,x_3) \mapsto (b_1x_2/(x_1x_3),b_2/x_3,b_3/x_2)
\end{eqnarray*}
where $b_1,b_2,b_3\in K\backslash \{0\}$.

Define $y_1=x_1$, $y_2=x_2/\sqrt{b_2}$, $y_3=x_3/\sqrt{b_2}$.
Use the fact that $\sigma\tau=\tau\sigma$.
We find that
\begin{eqnarray*}
\sigma &:& (y_1,y_2,y_3) \mapsto (1/y_1,y_3,y_2), \\
\tau &:& (y_1,y_2,y_3)\mapsto (\epsilon y_2/(y_1y_3),1/y_3,1/y_2)
\end{eqnarray*}
where $\epsilon=\pm 1$.

If $\epsilon=1$, apply Theorem \ref{t1.2}.
Thus $K(y_1,y_2,y_3)^G$ is rational.

It remains to consider the case that $\epsilon=-1$.

Define $t=y_2/y_3$. Then $K(y_1,y_2,y_3)=K(t,y_1,y_2)$ and
\begin{eqnarray*}
\sigma &:& (t,y_1,y_2) \mapsto (1/t,1/y_1,y_2/t), \\
\tau &:& (t,y_1,y_2)\mapsto (t,-t/y_1,t/y_2).
\end{eqnarray*}

Define
\[
u=\frac{y_1+\frac{t}{y_1}}{y_1y_2+\frac{t^2}{y_1y_2}},~~~
v=\frac{y_2-\frac{t}{y_2}}{y_1y_2+\frac{t^2}{y_1y_2}},~~~ w=u/v.
\]

By Theorem \ref{t2.3} we find that $K(t,y_1,y_2)^{\langle\tau\rangle}=K(t,u,w)$ and
\[
\sigma:t\mapsto 1/t,~w\mapsto w,~u\mapsto A/u
\]
where $A=w^2/(1+w^2)$ (the computation is similar to Case 4).

Define $s=[(1-t)/(1+t)][u-(A/u)]$. Then $K(t,u,w)=K(s,u,w)$ and
$\sigma(s)=s$. Thus we may reduce the question to the rationality
of $K(w)(u)^{\langle\sigma\rangle}$, which is easy by L\"uroth's
Theorem. Hence the result.
\end{idef}

\bigskip
\begin{idef}{Case 9.} $G=W_{13}(174)$.

The action of $G={\langle\sigma,\tau\rangle}$ is given by
\begin{eqnarray*}
\sigma &:& (x_1,x_2,x_3) \mapsto (1/x_1,x_3,x_2), \\
\tau &:& (x_1,x_2,x_3) \mapsto (a_1x_1x_3/x_2,a_2x_3,a_3x_2)
\end{eqnarray*}
where $a_1,a_2,a_3\in K\backslash \{0\}$.

Define $x_0=x_2/x_3$. Then $K(x_1,x_2,x_3)=K(x_0,x_1,x_2)$ and
\begin{eqnarray*}
\sigma &:& (x_0,x_1,x_2) \mapsto (1/x_0,1/x_1,x_2/x_0), \\
\tau &:& (x_0,x_1,x_2) \mapsto
(a_2/(a_3x_0),a_1x_1/x_0,a_2x_2/x_0)
\end{eqnarray*}

Apply Theorem \ref{t2.2}. It suffices to show that $K(x_0,x_1)^G$
is rational. But this follows from Theorem \ref{t1.1}.
\end{idef}

\bigskip
\begin{idef}{Case 10.} $G=W_{14}(174)$.

The action of $G={\langle\sigma,\tau\rangle}$ is given by
\begin{eqnarray*}
\sigma &:& (x_1,x_2,x_3) \mapsto (1/x_1,x_3,x_2), \\
\tau &:& (x_1,x_2,x_3) \mapsto (a_1/x_1,a_2x_1/x_3,a_3/(x_1x_2))
\end{eqnarray*}
where $a_1,a_2,a_3\in K\backslash \{0\}$.

Since $\tau^2=1$, we find that $a_3=a_1a_2$.
From $\sigma\tau(x_2)=\tau\sigma(x_2)$, we find that $a_2=a_3$.
Thus $a_1=1$. In summary, we have that
\begin{eqnarray*}
\sigma &:& (x_1,x_2,x_3) \mapsto (1/x_1,x_3,x_2), \\
\tau &:& (x_1,x_2,x_3) \mapsto (1/x_1,ax_1/x_3,a/(x_1x_2))
\end{eqnarray*}
where $a\in K\backslash \{0\}$.

Define $y_i=x_i/\sqrt{a}$ for $i=2,3$. Then
$\sigma:y_2\leftrightarrow y_3$, $\tau:y_2\mapsto x_1/y_3$,
$y_3\mapsto 1/(x_1y_2)$. Thus $G$ acts on $K(x_1,y_2,y_3)$ by
purely monomial $K$-automorphisms. Hence we may apply Theorem
\ref{t1.2}.
\end{idef}

\bigskip
\begin{idef}{Case 11.} $G=W_{15}(174)$.

The action of $G={\langle\sigma,\tau\rangle}$ is given by
\begin{eqnarray*}
\sigma &:& (x_1,x_2,x_3) \mapsto (1/x_1,x_3,x_2), \\
\tau &:& (x_1,x_2,x_3) \mapsto (a_1x_1,a_2x_3/x_1,a_3x_1x_2).
\end{eqnarray*}

Define $x_0=x_2/x_3$. The proof is almost the same as that of Case
9. $G=W_{13}(174)$.
\end{idef}
\end{proof}

\begin{theorem} \label{t3.5}
Let $G=W_i(187)$ for $1\le i\le 2$.
Then $K(x_1,x_2,x_3)^G$ is rational over $K$.
\end{theorem}

\begin{proof}
\begin{idef}{Case 1.} $G=W_1(187)$.

The action of $G={\langle\sigma,\tau\rangle}$ is given by
\begin{eqnarray*}
\sigma &:& (x_1,x_2,x_3) \mapsto \left( (\sqrt{-1})^i x_1, x_3,b/x_2 \right), \\
\tau &:& (x_1,x_2,x_3) \mapsto (c_1/x_1,c_2/x_2,c_3/x_3)
\end{eqnarray*}
where $b,c_1,c_2,c_3\in K\backslash\{0\}$ and $i\in \{0,1,2,3\}$.

Since $\sigma\tau=\tau\sigma$, we find that $i=0$ or 2, $c_2=c_3$, $(b/c_2)^2=1$.
In summary, we have
\begin{eqnarray*}
\sigma &:& (x_1,x_2,x_3) \mapsto (\epsilon_1x_1,x_3,a_2\epsilon_2/x_2), \\
\tau &:& (x_1,x_2,x_3) \mapsto (a_1/x_1,a_2/x_2,a_2/x_3)
\end{eqnarray*}
where $a_1,a_2\in K\backslash\{0\}$ and $\epsilon_1,\epsilon_2\in
\{1,-1\}$.

Define $y_1=x_1/\sqrt{a_1}$, $y_2=x_2/\sqrt{a_2}$, $y_3=x_3/\sqrt{a_2}$.
We get that
\begin{eqnarray*}
\sigma &:& (y_1,y_2,y_3) \mapsto (\epsilon_1y_1,y_3,\epsilon_2/y_2), \\
\tau &:& (y_1,y_2,y_3)\mapsto (1/y_1,1/y_2,1/y_3)
\end{eqnarray*}
where $\epsilon_1,\epsilon_2\in\{1,-1\}$.

If $\epsilon_1=1$, define $z_1=(1-y_1)/(1+y_1)$.
Then $\sigma(z_1)=z_1$, $\tau(z_1)=-z_1$.
Apply Theorem \ref{t2.1} and reduce the question to $K(y_2,y_3)^G$.

If $\epsilon_1=-1$ and $\epsilon_2=1$, define $z_2=(1-y_2)/(1+y_2)$, $z_3=(1-y_3)/(1+y_3)$.
Then $\sigma(z_2)=z_3$, $\sigma(z_3)=-z_2$, $\tau(z_2)=-z_2$, $\tau(z_3)=-z_3$.
Apply Theorem \ref{t2.1} and reduce the question to $K(y_1)^G$.

It remains to consider the situation when $\epsilon_1=\epsilon_2=-1$.

Note that $\sigma^2(y_1)=y_1$, $\sigma^2(y_2)=-1/y_2$, $\sigma^2(y_3)=-1/y_3$.
Define
\[
u=\frac{y_2+\frac{1}{y_2}}{y_2y_3-\frac{1}{y_2y_3}},~~~
v=\frac{y_3+\frac{1}{y_3}}{y_2y_3-\frac{1}{y_2y_3}},~~~
w=u/v.
\]

By Theorem \ref{t2.3} we find that $K(y_1,y_2,y_3)^{\langle\sigma^2\rangle}=K(y_1,u,v)=K(y_1,w,v)$ and
\begin{eqnarray*}
\sigma &:& (y_1,w,v)\mapsto (-y_1,-1/w,A/v) \\
\tau &:& (y_1,w,v)\mapsto (1/y_1,w,-v)
\end{eqnarray*}
where $A=w/(w^2-1)$.

Define $z_1=y_1+(1/y_1)$, $z_2=v[y_1-(1/y_1)]$.
Then $K(y_1,w,v)^{\langle\tau\rangle}=K(z_1,z_2,w)$ and
\[
\sigma: z_1\mapsto -z_1,~ w\mapsto -1/w,~ z_2\mapsto B/z_2
\]
where $B=(z_1^2-4)/[(1/w)-w]$.

Define $t=z_1[w+(1/w)]$, $z_3=z_2[(1/w)-w]/(z_1-2)$.
Then $K(z_1,z_2,w)=K(t,w,z_3)$ and
\[
\sigma: t\mapsto t,~ w\mapsto -1/w,~ z_3\mapsto [w-(1/w)]/z_3.
\]

Thus $K(t,w,z_3)^{\langle\sigma\rangle}$ is rational over $K(t)$ by Theorem \ref{t2.4}
(instead of Theorem \ref{t2.3}).
\end{idef}

\bigskip
\begin{idef}{Case 2.} $G=W_2(187)$.

The action of $G={\langle\sigma,\tau\rangle}$ is given by
\begin{eqnarray*}
\sigma &:& (x_1,x_2,x_3) \mapsto (a_1x_1,x_3,a_2x_1/x_2), \\
\tau &:& (x_1,x_2,x_3) \mapsto (b_1/x_1,b_2/x_2,b_3/x_3)
\end{eqnarray*}
where $a_i,b_j\in K\backslash\{0\}$.

Define $y_1=x_2$, $y_2=x_3$, $y_3=a_2x_1/x_2$.
Using the relation $\sigma^4=1$, $\sigma\tau=\tau\sigma$, we get
\begin{eqnarray*}
\sigma &:& (y_1,y_2,y_3) \mapsto (y_2,y_3,\epsilon y_1y_3/y_2), \\
\tau &:& (y_1,y_2,y_3)\mapsto (a/y_1,a/y_2,a/y_3)
\end{eqnarray*}
where $\epsilon \in\{1,-1\}$ and $a\in K\backslash\{0\}$.

Define $z_i=y_i/\sqrt{a}$ for $1\le i\le 3$. Then we have
\begin{eqnarray*}
\sigma &:& (z_1,z_2,z_3)\mapsto (z_2,z_3,\epsilon z_1z_3/z_2), \\
\tau &:& (z_1,z_2,z_3)\mapsto (1/z_1,1/z_2,1/z_3).
\end{eqnarray*}
where $\epsilon\in \{1,-1\}$.

If $\epsilon=1$, the action of $G$ is a purely monomial action.
Apply Theorem \ref{t1.2}.

It remains to consider the situation when $\epsilon=-1$.

Define $t=z_1z_3$. Then $K(z_1,z_2,z_3)=K(t,z_1,z_2)$ and
\begin{eqnarray*}
\sigma &:& (t,z_1,z_2)\mapsto (-t,z_2,t/z_1), \\
\tau &:& (t,z_1,z_2)\mapsto (1/t,1/z_1,1/z_2).
\end{eqnarray*}

Note that $\sigma^2(t)=t$, $\sigma^2(z_1)=t/z_1$, $\sigma^2(z_2)=-t/z_2$. Define
\[
u=\frac{z_1-\frac{t}{z_1}}{z_1z_2+\frac{t^2}{z_1z_2}},~~~
v=\frac{z_2+\frac{t}{z_2}}{z_1z_2+\frac{t^2}{z_1z_2}},~~~
w=u/v.
\]

By Theorem \ref{t2.3} we find that $K(t,z_1,z_2)^{\langle\sigma^2\rangle}=K(t,u,v)=K(t,w,v)$ and
\begin{eqnarray*}
\sigma &:& (t,w,v) \mapsto (-t,-1/w,A/v) \\
\tau &:& (t,w,v) \mapsto (1/t,-w,tv)
\end{eqnarray*}
where $A=-w/[t(w^2+1)]$.

Define $u_1=t+(1/t)$, $u_2=w[t-(1/t)]$, $u_3=v(1+t)$.
Then $K(t,w,v)^{\langle\tau\rangle}=K(u_1,u_2$, $u_3)$ and
\[
\sigma: u_1\mapsto -u_1,~u_2\mapsto (u_1^2-4)/u_2,~u_3\mapsto B/u_3
\]
where $B=A(1-t^2)=w(t^2-1)/[t(w^2+1)]=u_2(u_1^2-4)/[u_2^2+(u_1^2-4)]$.

Define $v_1=u_1$, $v_2=u_2/(u_1-2)$, $v_3=u_3$. Then
\[
\sigma: v_1\mapsto -v_1,~ v_2\mapsto -1/v_2,~ v_3\mapsto B/v_3
\]
where $B=(v_1+2)v_2/[v_2^2+((v_1+2)/(v_1-2))]$.

Define $v_4=(\sqrt{-1}-v_2)/(\sqrt{-1}+v_2)$.
Then $\sigma(v_4)=-v_4$.
It is not difficult to verify that
\[
B=(v_1^2-4)(v_4^2-1)/[4\sqrt{-1}(v_4^2+v_1v_4+1)].
\]

Define $s=v_1^2$, $\alpha=v_1$, $x=v_1v_4$, $y=e\cdot v_3(v_4^2+v_1v_4+1)/[(v_1-2)(v_4-1)]$ where $e\in K$
satisfying $e^2=4\sqrt{-1}$.
Note that $K(u_1,u_2,u_3)=K(\alpha,x,y)$ and we find that
\begin{equation}
\sigma: \alpha\mapsto -\alpha,~x\mapsto x,~y\mapsto f(x)/y \label{eq3.2}
\end{equation}
where $f(x)=\frac{1}{s}\{[x+(s/2)]^2-[(s^2-4s)/4]\}$.

We will apply Theorem \ref{t2.5} (note that Formula \eqref{eq3.2} looks very similar to Formula \eqref{eq3.1}
in the proof of Case 2 in Theorem \ref{t3.4}).
In the notation of Theorem \ref{t2.5}, $F=K(s)$ and $E=F(\alpha)$.
It remains to show that the Hilbert symbol $(s,1/s)_2$ is trivial in $K(s,\sqrt{s-4})$.
Since $-1\in K^2$, we find that $(s,1/s)_2=(s,s)_2=(s,-s)_2=0$.
Hence $K(\alpha,x,y)^{\langle\sigma\rangle}=K(s)(\alpha,x,y)^{\langle\sigma\rangle}$ is rational over $K(s)$.
\end{idef}
\end{proof}

\begin{theorem} \label{t3.6}
Let $G=W_i(187)$ for $3\le i \le 6$.
Then $K(x_1,x_2,x_3)^G$ is rational over $K$.
\end{theorem}

\begin{proof}
\begin{idef}{Case 1.} $G=W_3(187)$.

The action of $G=\langle\sigma,\tau,\lambda\rangle$ is given by
\begin{eqnarray*}
\sigma &:& (x_1,x_2,x_3)\mapsto (\epsilon_1x_1,1/x_2,1/x_3), \\
\tau &:& (x_1,x_2,x_3)\mapsto (1/x_1,\epsilon_2/x_2,\epsilon_3/x_3), \\
\lambda &:& (x_1,x_2,x_3)\mapsto (a_1/x_1,a_2/x_2,a_3/x_3)
\end{eqnarray*}
where $\epsilon_1,\epsilon_2,\epsilon_3\in \{1,-1\}$, $a_1,a_2,a_3\in K\backslash\{0\}$
(for the actions of $\sigma$ and $\tau$,
compare with the actions of $\sigma$ and $\tau$ for $G=W_6(174)$ in the proof of Theorem \ref{t3.4}).

Since $\sigma\lambda=\lambda\sigma$ and $\tau\lambda=\lambda\tau$, we find that $a_1^2=a_2^2=a_3^2=1$.

We will consider only the case $\epsilon_1=-1$, because the case
$\epsilon_1=1$ is easier and can be proved similarly.

Define $y_1=x_1$ and
\[
y_i=\frac{1-x_i}{1+x_i}
\]
for $i=2,3$. Thus we get
\begin{eqnarray*}
\sigma &:& (y_1,y_2,y_3)\mapsto (-y_1,-y_2,-y_3), \\
\tau &:& (y_1,y_2,y_3)\mapsto (y_1^{-1},-y_2^{\epsilon_2},-y_3^{\epsilon_3}), \\
\lambda &:& (y_1,y_2,y_3)\mapsto
(a_1y_1^{-1},-y_2^{a_2},-y_3^{a_3})
\end{eqnarray*}
where $\epsilon_2,\epsilon_3,a_1,a_2,a_3\in\{1,-1\}$.

Define $z_1=y_1/y_2$, $z_2=y_2/y_3$, $z_3=y_3^2$.
Then $K(y_1,y_2,y_3)^{\langle\sigma\rangle}=K(z_1,z_2,z_3)$.
It is not difficult to see that anyone of $\tau(z_i)$, $\lambda(z_j)$ where $1\le i,j\le 3$ is of the form
\[
\epsilon z_1^{b_1} z_2^{b_2}z_3^{b_3}
\]
where $\epsilon\in \{1,-1\}$ and $b_1,b_2,b_3\in \bm{Z}$. Hence
$\langle\tau,\lambda\rangle$ acts on $K(z_1,z_2,z_3)$ by monomial
$K$-automorphisms. Since $\langle\tau,\lambda\rangle \simeq
C_2\times C_2$, we find that
$K(z_1,z_2,z_3)^{\langle\tau,\lambda\rangle}$ is rational over $K$
by Theorem \ref{t3.4}. Done.
\end{idef}

\bigskip
\begin{idef}{Case 2.} $G=W_4(187)$.

The action of $G=\langle\sigma,\tau,\lambda\rangle$ is given by
\begin{eqnarray*}
\sigma &:& (x_1,x_2,x_3)\mapsto (a_1x_1,a_2/x_2,a_3/x_3), \\
\tau &:& (x_1,x_2,x_3)\mapsto (b_1/x_1,b_2/x_3,b_3/x_2), \\
\lambda &:& (x_1,x_2,x_3)\mapsto (c_1/x_1,c_2/x_2,c_3/x_3)
\end{eqnarray*}
where $a_i,b_j,c_k \in K\backslash\{0\}$.

Define $y_i=x_i/\sqrt{c_i}$ for $1\le i\le 3$.
It is easy to see that
\begin{eqnarray*}
\sigma &:& (y_1,y_2,y_3)\mapsto (\epsilon_1y_1,\epsilon_2/y_2,\epsilon_3/y_3), \\
\tau &:& (y_1,y_2,y_3)\mapsto (e_1/y_1,e_2/y_3,e_3/y_2), \\
\lambda &:& (y_1,y_2,y_3)\mapsto (1/y_1,1/y_2,1/y_3)
\end{eqnarray*}
where $\epsilon_1,\epsilon_2,\epsilon_3,e_1,e_2,e_3\in\{1,-1\}$ and $e_2=e_3$.

Define $z_i=(1-y_i)/(1+y_i)$ for $1\le i\le 3$. Consider
$K(z_1,z_2,z_3)^{\langle\lambda\rangle}$ and then proceed as in
the above Case 1. The details are left to the reader.
\end{idef}

\bigskip
\begin{idef}{Case 3.} $G=W_5(187)$.

The action of $G=\langle\sigma,\tau,\lambda\rangle$ is given by
\begin{eqnarray*}
\sigma &:& (x_1,x_2,x_3)\mapsto (1/x_1,x_3,x_2), \\
\tau &:& (x_1,x_2,x_3)\mapsto (\epsilon x_2/(x_1x_3),1/x_3,1/x_2), \\
\lambda &:& (x_1,x_2,x_3)\mapsto (a_1/x_1,a_2/x_2,a_3/x_3)
\end{eqnarray*}
where $\epsilon\in \{1,-1\}$, $a_1,a_2,a_3\in K\backslash\{0\}$
(for the actions of $\sigma$ and $\tau$, compare with those for
$G=W_{12}(174)$ in the proof of Theorem \ref{t3.4}).

Since $G\simeq C_2\times C_2\times C_2$, it is easy to see that
$a_2=a_3$ and $a_1,a_2\in \{1,-1\}$.

Define $y_1=x_1$, $y_2=x_2$, $y_3=x_2/x_3$. We find that
\begin{eqnarray*}
\sigma &:& (y_1,y_2,y_3)\mapsto (1/y_1,y_2/y_3,1/y_3), \\
\tau &:& (y_1,y_2,y_3)\mapsto (\epsilon y_3/y_1,y_3/y_2,y_3), \\
\lambda &:& (y_1,y_2,y_3)\mapsto (\epsilon_1/y_1,\epsilon_2/y_2,1/y_3)
\end{eqnarray*}
where $\epsilon,\epsilon_1,\epsilon_2\in\{1,-1\}$.

Define
\[
u=\frac{y_1-\frac{\epsilon y_3}{y_1}}{y_1y_2-\frac{\epsilon y_3^2}{y_1y_2}},~~~
v=\frac{y_2-\frac{y_3}{y_2}}{y_1y_2-\frac{\epsilon y_3^2}{y_1y_2}},~~~ w=u/v.
\]

By Theorem \ref{t2.3}, we get $K(y_1,y_2,y_3)^{\langle\tau\rangle}=K(u,v,y_3)$ and
\begin{eqnarray*}
\sigma &:& (u,w,y_3)\mapsto (A/u,-\epsilon w,1/y_3), \\
\lambda &:& (u,w,y_3)\mapsto (\epsilon_2y_3u,\epsilon\epsilon_1\epsilon_2w,1/y_3)
\end{eqnarray*}
where $A=w^2/(w^2-\epsilon)$.

Define $z_1=w$, $z_2=(1-y_3)/(1+y_3)$, $z_3=2uy_3/(1+y_3)$.
Then $K(u,w,y_3)=K(z_1,z_2,z_3)$ and
\begin{eqnarray*}
\sigma &:& (z_1,z_2,z_3)\mapsto (-\epsilon z_1,-z_2,B/z_3), \\
\lambda &:& (z_1,z_2,z_3)\mapsto (\epsilon\epsilon_1\epsilon_2z_1,-z_2,\epsilon_2z_3)
\end{eqnarray*}
where $B=A(1-z_2^2)=z_1^2(1-z_2^2)/(z_1^2-\epsilon)$.

\bigskip
\begin{idef}{Case 3.1.} $\epsilon_2=-1$ and $\epsilon\epsilon_1\epsilon_2=-1$.

Define $u_1=z_1z_2$, $u_2=z_2^2$, $u_3=z_2z_3$.
Then $K(z_1,z_2,z_3)^{\langle\lambda\rangle}=K(u_1,u_2,u_3)$ and
\[
\sigma: (u_1,u_2,u_3) \mapsto (\epsilon u_1,u_2,C/u_3)
\]
where $C=u_1^2u_2(1-u_2)/(\epsilon-u_1^2)$.

If $\epsilon=1$, then $K(u_1,u_2,u_3)^{\langle\sigma\rangle}=K(u_1,u_2,u_3+(C/u_3))$.

If $\epsilon=-1$, define $u_4=u_3(1-\sqrt{-1}u_1)/u_1$. Then we get that
\[
\sigma : (u_1,u_2,u_4)\mapsto (-u_1,u_2,u_2(1-u_2)/u_4).
\]

Thus $\sigma$ acts on $K(u_2)(u_1,u_4)$ as a monomial automorphism
over $K(u_2)$. Hence
$K(u_1,u_2,u_4)^{\langle\sigma\rangle}=K(u_2)(u_1,u_4)^{\langle\sigma\rangle}$
is rational over $K(u_2)$.
\end{idef}

\bigskip
\begin{idef}{Case 3.2.} $\epsilon_2=-1$ and $\epsilon\epsilon_1\epsilon_2=1$.

Define $u_1=z_1$, $u_2=z_2^2$, $u_3=z_2z_3$. The details are left
to the reader.
\end{idef}

\bigskip
\begin{idef}{Case 3.3.} $\epsilon_2=1$.

This case is easier and the proof is omitted.
\end{idef}
\end{idef}

\bigskip
\begin{idef}{Case 4.} $G=W_6(187)$.

The action of $G=\langle\sigma,\tau,\lambda\rangle$ is given by
\begin{eqnarray*}
\sigma &:& (x_1,x_2,x_3)\mapsto (1/x_1,x_3,x_2), \\
\tau &:& (x_1,x_2,x_3)\mapsto (1/x_1,x_1/x_3,1/(x_1x_2)), \\
\lambda &:& (x_1,x_2,x_3)\mapsto (a_1/x_1,a_2/x_2,a_3/x_3)
\end{eqnarray*}
where $a_1,a_2,a_3\in K\backslash\{0\}$
(for the action of $\sigma$ and $\tau$,
see $G=W_{14}(174)$ in the proof of Theorem \ref{t3.4}).

Since $G\simeq C_2\times C_2\times C_2$, it follows that $a_1^2=1$, $a_2=a_3$, and $a_2^2=a_1$.

Define $y_1=x_1x_2/x_3$, $y_2=x_2$, $y_3=x_3$. Then we have
\begin{eqnarray*}
\sigma &:& (y_1,y_2,y_3)\mapsto (1/y_1,y_3,y_2), \\
\tau &:& (y_1,y_2,y_3)\mapsto (y_1,y_1/y_2,1/(y_1y_3)), \\
\lambda &:& (y_1,y_2,y_3)\mapsto (\epsilon/y_1,a/y_2,a/y_3)
\end{eqnarray*}
where $\epsilon\in\{1,-1\}$ and $a^2=\epsilon$.

Define
\[
u=\frac{y_2-\frac{y_1}{y_2}}{y_2y_3-\frac{1}{y_2y_3}},~~~
v=\frac{y_3-\frac{1}{y_1y_3}}{y_2y_3-\frac{1}{y_2y_3}}.
\]

By Theorem \ref{t2.3}, we find that $K(y_1,y_2,y_3)^{\langle\tau\rangle}=K(y_1,u,v)$ and
\begin{eqnarray*}
\sigma &:& (y_1,u,v)\mapsto (1/y_1,v,u), \\
\lambda &:& (y_1,u,v)\mapsto (\epsilon/y_1,u/(ay_1),y_1v/a).
\end{eqnarray*}

Thus $\langle\sigma,\lambda\rangle$ acts on $K(y_1,u,v)$ by
monomial $K$-automorphisms. By Theorem \ref{t3.4},
$K(y_1,u,v)^{\langle\sigma,\lambda\rangle}$ is rational over $K$.
\end{idef}
\end{proof}

\begin{theorem}\label{t3.7}
Let $G=W_i(187)$ for $7\le i \le 14$.
Then $K(x_1,x_2,x_3)^G$ is rational over $K$.
\end{theorem}

\begin{proof}
\begin{idef}{Case 1.} $G=W_7(187)$.

The action of $G=\langle\sigma,\tau\rangle$ is given by
\begin{eqnarray*}
\sigma &:& (x_1,x_2,x_3) \mapsto \left( (\sqrt{-1})^j x_1,x_3,b/x_2 \right), \\
\tau &:& (x_1,x_2,x_3) \mapsto (b_1/x_1,b_2x_3,b_3x_2)
\end{eqnarray*}
where $b,b_1,b_2,b_3\in K\backslash\{0\}$ and $j\in \{0,1,2,3\}$
(for the action of $\sigma$, see Case 1. $G=W_1(187)$ in the proof of Theorem \ref{t3.5}).

Define $y_1=x_1/\sqrt{b_1}$, $y_2=x_2/\sqrt{b}$, $y_3=b_2x_3/\sqrt{b}$.
Then we have
\begin{eqnarray*}
\sigma &:& (y_1,y_2,y_3) \mapsto \left( (\sqrt{-1})^j y_1,y_3/a,a/y_2 \right) \\
\tau &:& (y_1,y_2,y_3)\mapsto (1/y_1,y_3,y_2)
\end{eqnarray*}
for some $a\in K\backslash\{0\}$ and $j\in \{0,1,2,3\}$.

Since $\tau\sigma\tau^{-1}=\sigma^{-1}$, we find that $a^2=1$.

For $2\le i\le 3$, define $z_i=y_i/\sqrt{-1}$ if $a=-1$, and define $z_i=y_i$ if $a=1$.
Thus we get $K(x_1,x_2,x_3)=K(y_1,z_2,z_3)$ and
\begin{eqnarray*}
\sigma &:& (y_1,z_2,z_3)\mapsto \left( (\sqrt{-1})^j y_1,\epsilon z_3,1/z_2 \right), \\
\tau &:& (y_1,z_2,z_3)\mapsto (1/y_1,z_3,z_2)
\end{eqnarray*}
where $\epsilon\in\{1,-1\}$ and $j\in \{0,1,2,3\}$.

\bigskip
\begin{idef}{Case 1.1.} $\epsilon=1$.

Define $u_i=(1-z_i)/(1+z_i)$ for $2\le i\le 3$. We get
\begin{eqnarray*}
\sigma &:& (y_1,u_2,u_3)\mapsto \left( (\sqrt{-1})^j y_1,u_3,-u_2 \right), \\
\tau &:& (y_1,u_2,u_3)\mapsto (1/y_1,u_3,u_2).
\end{eqnarray*}

If $j\in \{1,3\}$, apply Theorem \ref{t2.1} and reduce the question to $K(y_1)^G$.

If $j\in \{0,2\}$, define $w_2=u_2/u_3$, $w_3=u_2^2$. Then
$K(y_1,u_2,u_3)^{\langle\sigma^2\rangle}=K(y_1,w_2,w_3)$. Note
that $\langle\sigma,\tau\rangle$ acts on $K(y_1,w_2,w_3)$ by
monomial $K$-automorphisms and $\langle\sigma,\tau\rangle \simeq
C_2\times C_2$ as a group of automorphisms on $K(y_1,w_2,w_3)$.
Thus we may apply Theorem \ref{t3.4}.
\end{idef}

\bigskip
\begin{idef}{Case 1.2.} $\epsilon=-1$ and $j=0$.

Define $z_1=(1-y_1)/(1+y_1)$.
Then $\sigma(z_1)=z_1$, $\tau(z_1)=-z_1$.
Apply Theorem \ref{t2.1} and it suffices to consider the rationality of $K(z_2,z_3)^G$.
\end{idef}

\bigskip
\begin{idef}{Case 1.3.} $\epsilon=-1$ and $j=2$.

Define
\[
u=\frac{z_2+\frac{1}{z_2}}{z_2z_3-\frac{1}{z_2z_3}},~~~
v=\frac{z_3+\frac{1}{z_3}}{z_2z_3-\frac{1}{z_2z_3}}.
\]

By Theorem \ref{t2.3}, $K(y_1,z_2,z_3)^{\langle\sigma^2\rangle}=K(y_1,u,v)$ and
\begin{eqnarray*}
\sigma &:& (y_1,u,v)\mapsto (-y_1,-v/(u^2-v^2),u/(u^2-v^2)) \\
\tau &:& (y_1,u,v)\mapsto (1/y_1,v,u).
\end{eqnarray*}

Define $w_1=u+v$, $w_2=u-v$.
Then we find that
\begin{eqnarray*}
\sigma &:& (y_1,w_2,w_3)\mapsto (-y_1,1/w_1,-1/w_2), \\
\tau &:& (y_1,w_2,w_3)\mapsto (1/y_1,w_1,-w_2).
\end{eqnarray*}

We get a monomial group action.
Thus $K(y_1,w_2,w_3)^{\langle\sigma,\tau\rangle}$ is rational by Theorem \ref{t3.4}.
\end{idef}

\bigskip
\begin{idef}{Case 1.4.} $\epsilon=-1$ and $j\in\{1,3\}$. Write $\eta=(\sqrt{-1})^j$.

Define $u_2=z_2/\eta$, $u_3=-z_3/\eta$. We get
\begin{eqnarray*}
\sigma &:& (y_1,u_2,u_3)\mapsto (\eta y_1,u_3,1/u_2), \\
\tau &:& (y_1,u_2,u_3)\mapsto (1/y_1,-u_3,-u_2).
\end{eqnarray*}

Define $v_2=(1-u_2)/(1+u_2)$, $v_3=(1-u_3)/(1+u_3)$.
We find that $K(y_1,z_2,z_3)=K(y_1,v_2,v_3)$ and
\begin{eqnarray*}
\sigma &:& (y_1,v_2,v_3)\mapsto (\eta y_1,v_3,-v_2), \\
\tau &:& (y_1,v_2,v_3)\mapsto (1/y_1,1/v_3,1/v_2), \\
\sigma^2 &:& (y_1,v_2,v_3) \mapsto (-y_1,-v_2,-v_3).
\end{eqnarray*}

Define $w_1=y_1v_2$, $w_2=v_2/v_3$, $w_3=v_2^2$.
Then $K(y_1,v_2,v_3)^{\langle\sigma^2\rangle}=K(w_1,w_2,w_3)$ and $\langle\sigma,\tau\rangle$ acts on
$K(w_1,w_2,w_3)$ by monomial $K$-automorphisms because
\begin{eqnarray*}
\sigma &:& (w_1,w_2,w_3) \mapsto (\eta w_1/w_2,-1/w_2,w_3/w_2^2), \\
\tau &:& (w_1,w_2,w_3) \mapsto (-w_3/w_1,w_2,w_3).
\end{eqnarray*}

Thus we may apply Theorem \ref{t3.4}.
\end{idef}
\end{idef}

\bigskip
\begin{idef}{Case 2.} $G=W_8(187)$.

The action of $G=\langle\sigma,\tau\rangle$ is given by
\begin{eqnarray*}
\sigma &:& (x_1,x_2,x_3)\mapsto (a_1x_1,a_2x_3,a_3/x_2), \\
\tau &:& (x_1,x_2,x_3)\mapsto (b_1x_1,b_2/x_3,b_3/x_2)
\end{eqnarray*}
where $a_i,b_j\in K\backslash \{0\}$.

Apply Theorem \ref{t2.1} and reduce the question to $K(x_2,x_3)^G$.
\end{idef}

\bigskip
\begin{idef}{Case 3.} $G=W_9(187)$.

The action of $G=\langle\sigma,\tau\rangle$ is given by
\begin{eqnarray*}
\sigma &:& (x_1,x_2,x_3)\mapsto (1/x_1,b_2/x_3,b_3x_2), \\
\tau &:& (x_1,x_2,x_3)\mapsto (c_1/x_1,c_2x_3,c_3x_2)
\end{eqnarray*}
where $b_2,b_3,c_1,c_2,c_3\in K\backslash \{0\}$.

Define $y_1=x_1$, $y_2=x_2$, $y_3=b_2/x_3$. Use the fact that
$G\simeq D_4$. We find that
\begin{eqnarray*}
\sigma &:& (y_1,y_2,y_3)\mapsto (1/y_1,y_3,a/y_2), \\
\tau &:& (y_1,y_2,y_3)\mapsto (\epsilon_1/y_1,a\epsilon_2/y_3,a\epsilon_2/y_2)
\end{eqnarray*}
where $a\in K\backslash \{0\}$ and $\epsilon_1,\epsilon_2\in \{1,-1\}$.

Define $z_i=y_i/\sqrt{a}$ for $2\le i \le 3$.
Then $K(x_1,x_2,x_3)=K(y_1,z_2,z_3)$ and
\begin{eqnarray*}
\sigma &:& (y_1,z_2,z_3) \mapsto (1/y_1,z_3,1/z_2), \\
\tau &:& (y_1,z_2,z_3)\mapsto (\epsilon_1/y_1,\epsilon_2/z_3,\epsilon_2/z_2)
\end{eqnarray*}
where $\epsilon_1,\epsilon_2\in \{1,-1\}$.

Define
\[
u=\frac{z_2-\frac{1}{z_2}}{z_2z_3-\frac{1}{z_2z_3}}, ~~~
v=\frac{z_3-\frac{1}{z_3}}{z_2z_3-\frac{1}{z_2z_3}}.
\]

By Theorem \ref{t2.3}, $K(y_1,z_2,z_3)^{\langle\sigma^2\rangle}=K(y_1,u,v)$ and
\begin{eqnarray*}
\sigma &:& (y_1,u,v) \mapsto (1/y_1,-v/(u^2-v^2),u/(u^2-v^2)), \\
\tau &:& (y_1,u,v)\mapsto
(\epsilon_1/y_1,\epsilon_2v,\epsilon_2u).
\end{eqnarray*}

Define $w_2=u+v$, $w_3=u-v$. Then
\begin{eqnarray*}
\sigma &:& (y_1,w_2,w_3) \mapsto (1/y_1,1/w_2,-1/w_3), \\
\tau &:& (y_1,w_2,w_3) \mapsto (\epsilon_1/y_1,\epsilon_2w_2,-\epsilon_2w_3).
\end{eqnarray*}

We may apply Theorem \ref{t3.4} to assert that $K(y_1,w_2,w_3)^{\langle\sigma,\tau\rangle}$ is rational over $K$.
\end{idef}

\bigskip
\begin{idef}{Case 4.} $G=W_{10}(187)$.

The action of $G=\langle\sigma,\tau\rangle$ is given by
\begin{eqnarray*}
\sigma &:& (x_1,x_2,x_3)\mapsto (1/x_1,a_2/x_3,a_3x_2), \\
\tau &:& (x_1,x_2,x_3)\mapsto (b_1x_1,b_2/x_3,b_3/x_2)
\end{eqnarray*}
where $a_2,a_3,b_1,b_2,b_3\in K\backslash \{0\}$.

Define $y_1=x_1$, $y_2=\sqrt{a_3}x_2/\sqrt{a_2}$, $y_3=(\sqrt{a_2}\cdot\sqrt{a_3})/x_3$.
Use the relations that $\tau^2=1$ and $\tau\sigma\tau^{-1}=\sigma^{-1}$. We get
\begin{eqnarray*}
\sigma &:& (y_1,y_2,y_3)\mapsto (1/y_1,y_3,1/y_2), \\
\tau &:& (y_1,y_2,y_3)\mapsto (\epsilon_1y_1,\epsilon_2y_3,\epsilon_2y_2)
\end{eqnarray*}
where $\epsilon_1,\epsilon_2\in \{1,-1\}$.

Define
\[
u=\frac{y_2-\frac{1}{y_2}}{y_2y_3-\frac{1}{y_2y_3}}, ~~~
v=\frac{y_3-\frac{1}{y_3}}{y_2y_3-\frac{1}{y_2y_3}}.
\]

By Theorem \ref{t2.3} we find that $K(y_1,y_2,y_3)^{\langle\sigma^2\rangle}=K(y_1,u,v)$ and
\begin{eqnarray*}
\sigma &:& (y_1,u,v) \mapsto (1/y_1,-v/(u^2-v^2),u/(u^2-v^2)), \\
\tau &:& (y_1,u,v)\mapsto (\epsilon_1y_1,\epsilon_2v,\epsilon_2u).
\end{eqnarray*}

Define $z_2=u+v$, $z_3=u-v$. Then $\langle\sigma,\tau\rangle$ acts on $K(y_1,u,v)=K(y_1,z_2,z_3)$ by
monomial $K$-automorphisms because
\begin{eqnarray*}
\sigma &:& (y_1,z_2,z_3) \mapsto (1/y_1,1/z_2,-1/z_3), \\
\tau &:& (y_1,z_2,z_3) \mapsto (\epsilon_1y_1,\epsilon_2z_2,-\epsilon_2z_3).
\end{eqnarray*}

Thus we may apply Theorem \ref{t3.4}.
\end{idef}

\bigskip
\begin{idef}{Case 5.} $G=W_{11}(187)$.

The action of $G=\langle\sigma,\tau\rangle$ is given by
\begin{eqnarray*}
\sigma &:& (x_1,x_2,x_3)\mapsto (a_1x_1,x_3,a_3x_1/x_2), \\
\tau &:& (x_1,x_2,x_3)\mapsto (1/x_1,b_2/x_3,b_3/x_2)
\end{eqnarray*}
where $a_1,a_3,b_2,b_3\in K\backslash \{0\}$.

Use the relations that $\tau^4=\tau^2=1$ and $\tau\sigma\tau^{-1}=\sigma^{-1}$. We find that
\begin{eqnarray*}
\sigma &:& (x_1,x_2,x_3)\mapsto (\epsilon x_1,x_3,a_3x_1/x_2), \\
\tau &:& (x_1,x_2,x_3)\mapsto (1/x_1,b/x_3,b/x_2)
\end{eqnarray*}
where $\epsilon\in\{1,-1\}$, $a_3,b\in K\backslash \{0\}$.

Define $y_1=x_1$, $y_2=x_2/\sqrt{b}$, $y_3=\sqrt{b}/x_3$ and use that relation
$\tau\sigma\tau^{-1}(y_3)=\sigma^{-1}(y_3)$. We find that
\begin{eqnarray*}
\sigma &:& (y_1,y_2,y_3)\mapsto (\epsilon y_1,1/y_3,ay_2/y_1), \\
\tau &:& (y_1,y_2,y_3)\mapsto (1/y_1,y_3,y_2)
\end{eqnarray*}
where $\epsilon\in \{1,-1\}$ and $a^2=\epsilon$.

Note that $\sigma^2: (y_1,y_2,y_3)\mapsto (y_1,y_1/(ay_2),1/(ay_1y_3))$.

Define
\[
u=\frac{y_2-\frac{y_1}{ay_2}}{y_2y_3-\frac{\epsilon}{y_2y_3}}, ~~~
v=\frac{y_3-\frac{1}{ay_1y_3}}{y_2y_3-\frac{\epsilon}{y_2y_3}}.
\]

By Theorem \ref{t2.3},
$K(y_1,y_2,y_3)^{\langle\sigma^2\rangle}=K(y_1,u,v)$ (note that
$a^2=\epsilon$) and
\begin{eqnarray*}
\sigma &:& (y_1,u,v) \mapsto (\epsilon y_1,-y_1^2v/(u^2-y_1^2v^2),u/(u^2-y_1^2v^2)) \\
\tau &:& (y_1,u,v)\mapsto (1/y_1,v,u).
\end{eqnarray*}

\bigskip
\begin{idef}{Case 5.1.} $\epsilon=-1$.

Define $z_2=u+y_1v$, $z_3=u-y_1v$. Then $K(y_1,u,v)=K(y_1,z_2,z_3)$ and
\begin{eqnarray*}
\sigma &:& (y_1,z_2,z_3) \mapsto (-y_1,-y_1/z_3,y_1/z_2), \\
\tau &:& (y_1,z_2,z_3) \mapsto (1/y_1,z_2/y_1,-z_3/y_1).
\end{eqnarray*}

Apply Theorem \ref{t3.4}.
\end{idef}

\bigskip
\begin{idef}{Case 5.2.} $\epsilon=1$.

Define $z_2$ and $z_3$ as in Case 1. The details are omitted.
\end{idef}
\end{idef}

\bigskip
\begin{idef}{Case 6.} $G=W_{12}(187)$.

The action of $G=\langle\sigma,\tau\rangle$ is given by
\begin{eqnarray*}
\sigma &:& (x_1,x_2,x_3)\mapsto (a_1x_1,a_2x_3,a_3x_1/x_2), \\
\tau &:& (x_1,x_2,x_3)\mapsto (\epsilon x_1,x_3,x_2)
\end{eqnarray*}
where $a_1,a_2,a_3\in K\backslash \{0\}$ and $\epsilon\in \{1,-1\}$.

Since $\tau^4=1$ and $\tau\sigma\tau^{-1}=\sigma^{-1}$, we find that $a_1=\epsilon$, $a_2^2=1$.

Define $y_1=x_1$, $y_2=x_2/\sqrt{a_3}$, $y_3=x_3/\sqrt{a_3}$. We get that
\begin{eqnarray*}
\sigma &:& (y_1,y_2,y_3)\mapsto (\epsilon_1 y_1,\epsilon_2y_3,y_1/y_2), \\
\tau &:& (y_1,y_2,y_3)\mapsto (\epsilon_1y_1,y_3,y_2)
\end{eqnarray*}
where $\epsilon_1,\epsilon_2\in \{1,-1\}$.

Define $z_1=y_1$, $z_2=y_2$, $z_3=1/y_3$. Then we have
\begin{eqnarray*}
\sigma &:& (z_1,z_2,z_3) \mapsto (\epsilon_1z_1,\epsilon_2/z_3,z_2/z_1), \\
\tau &:& (z_1,z_2,z_3) \mapsto (\epsilon_1z_1,1/z_3,1/z_2), \\
\sigma^2 &:& (z_1,z_2,z_3) \mapsto (z_1,\epsilon_2z_1/z_2,\epsilon_1\epsilon_2/(z_1z_3)).
\end{eqnarray*}

Define
\[
u=\frac{z_2-\frac{\epsilon_2z_1}{z_2}}{z_2z_3-\frac{\epsilon_1}{z_2z_3}}, ~~~
v=\frac{z_3-\frac{\epsilon_1\epsilon_2}{z_1z_3}}{z_2z_3-\frac{\epsilon_1}{z_2z_3}}.
\]

Then proceed by the same way as in the previous Case 5.
$G=W_{11}(187)$. Note that
$\sigma:(u,v)\mapsto(-\epsilon_1\epsilon_2
z_1^2v/(u^2-\epsilon_1z_1^2v^2)$,
$\epsilon_2u/(u^2-\epsilon_1z_1^2v^2))$. Thus, if $\epsilon_1=-1$,
we should take the change of variables $u\pm \sqrt{-1}z_1v$,
instead of the ``usual" change of variables $u\pm z_1v$. The
details are left to the reader.
\end{idef}

\bigskip
\begin{idef}{Case 7.} $G=W_{13}(187)$.

The action of $G=\langle\sigma,\tau\rangle$ is given by
\begin{eqnarray*}
\sigma &:& (x_1,x_2,x_3)\mapsto (1/x_1,1/x_3,a_3x_2/x_1), \\
\tau &:& (x_1,x_2,x_3)\mapsto (b_1/x_1,b_2/x_3,b_3/x_2)
\end{eqnarray*}
where $a_3,b_1,b_2,b_3\in K\backslash \{0\}$.

Use the relations $\tau^2=1$ and $\tau\sigma\tau^{-1}=\sigma^{-1}$. We find that $b_1=1$ and $b_2=b_3\in\{1$, $-1\}$.
Define $y_1=x_1$, $y_2=\sqrt{a_3}x_2$, $y_3=x_3/\sqrt{a_3}$. Then we get
\begin{eqnarray*}
\sigma &:& (y_1,y_2,y_3)\mapsto (1/y_1,1/y_3,y_2/y_1), \\
\tau &:& (y_1,y_2,y_3)\mapsto (1/y_1,\epsilon/y_3,\epsilon/y_2)
\end{eqnarray*}
where $\epsilon\in \{1,-1\}$.

If $\epsilon=1$, apply Theorem \ref{t1.2}.

If $\epsilon=-1$, note that $\sigma^2:(y_1,y_2,y_3)\mapsto (y_1,y_1/y_2,y_1/y_3)$.
Define
\begin{equation}
u=\frac{y_2-\frac{y_1}{y_2}}{y_2y_3-\frac{y_1^2}{y_2y_3}}, ~~~
v=\frac{y_3-\frac{y_1}{y_3}}{y_2y_3-\frac{y_1^2}{y_2y_3}}. \label{eq3.3}
\end{equation}

The proof is almost the same as the previous Case 5.
$G=W_{11}(187)$ by taking the change of variables $z_2=u+v$,
$z_3=u-v$ because $\sigma:(u,v)\mapsto (-v/(u^2-v^2),u/(u^2-v^2))$
and $\tau:(u,v)\mapsto (-y_1v,-y_1u)$.
\end{idef}

\bigskip
\begin{idef}{Case 8.} $G=W_{14}(187)$.

The action of $G=\langle\sigma,\tau\rangle$ is given by
\begin{eqnarray*}
\sigma &:& (x_1,x_2,x_3)\mapsto (1/x_1,1/x_3,a_3x_2/x_1), \\
\tau &:& (x_1,x_2,x_3)\mapsto (\epsilon x_1,bx_3,x_2/b)
\end{eqnarray*}
where $a_3,b\in K\backslash \{0\}$, $\epsilon\in\{1,-1\}$.

Define $y_1=x_1$, $y_2=\sqrt{a_3}x_2$, $y_3=x_3/\sqrt{a_3}$. Use
the relation $\tau\sigma\tau^{-1}=\sigma^{-1}$. We get
\begin{eqnarray*}
\sigma &:& (y_1,y_2,y_3)\mapsto (1/y_1,1/y_3,y_2/y_1), \\
\tau &:& (y_1,y_2,y_3)\mapsto (\epsilon y_1,ay_3,y_2/a)
\end{eqnarray*}
where $a^2=\epsilon\in \{1,-1\}$.

Define $u$ and $v$ by Formula \eqref{eq3.3} in the previous Case 7. $G=W_{13}(187)$.
Note that $\sigma:(y_1,u,v)\mapsto (1/y_1,-v/(u^2-v^2),u/(u^2-v^2))$,
$\tau:(y_1,u,v)\mapsto (\epsilon y_1,av,\epsilon au)$.

Discuss the situations $\epsilon=1$ and $\epsilon=-1$ separately.
The proof is almost the same as in the Case 7. $G=W_{13}(187)$.
\end{idef}
\end{proof}

\begin{theorem} \label{t3.8}
Let $G=W_1(194), W_2(195)$.
Then $K(x_1,x_2,x_3)^G$ is rational over $K$.
\end{theorem}

\begin{proof}
\begin{idef}{Case 1.} $G=W_1(194)$.

The action of $G=\langle\sigma,\tau,\lambda\rangle$ is given by
\begin{eqnarray*}
\sigma &:& (x_1,x_2,x_3)\mapsto \left( (\sqrt{-1})^j x_1,\epsilon x_3,1/x_2\right), \\
\tau &:& (x_1,x_2,x_3)\mapsto (1/x_1,x_3,x_2), \\
\lambda &:& (x_1,x_2,x_3) \mapsto (a_1/x_1,a_2/x_2,a_3/x_3)
\end{eqnarray*}
where $j\in \{0,1,2,3\}$, $\epsilon\in\{1,-1\}$, $a_1,a_2,a_3\in K\backslash \{0\}$
(for the actions of $\sigma$ and $\tau$, see Case 1. $G=W_7(187)$ in the proof of Theorem \ref{t3.7}).

Use the relations $\sigma\lambda=\lambda\sigma$ and $\tau\lambda=\lambda\tau$.
We find that $j=0$ or 2, $a_2=a_3\in \{1,-1\}$.
We write these actions as follows
\begin{eqnarray*}
\sigma &:& (x_1,x_2,x_3)\mapsto (\epsilon_1x_1,\epsilon_2 x_3,1/x_2), \\
\tau &:& (x_1,x_2,x_3)\mapsto (1/x_1,x_3,x_2), \\
\lambda &:& (x_1,x_2,x_3) \mapsto (a/x_1,\epsilon_3/x_2,\epsilon_3/x_3)
\end{eqnarray*}
where $\epsilon_1,\epsilon_2,\epsilon_3\in \{1,-1\}$ and $a\in K\backslash \{0\}$.

Define
\[
u=\frac{x_2-\frac{\epsilon_2}{x_2}}{x_2x_3-\frac{1}{x_2x_3}},~~~
v=\frac{x_3-\frac{\epsilon_2}{x_3}}{x_2x_3-\frac{1}{x_2x_3}}.
\]

By Theorem \ref{t2.3} we find that $K(x_1,x_2,x_3)^{\langle\sigma^2\rangle}=K(x_1,u,v)$ and
\begin{eqnarray*}
\sigma &:& (x_1,u,v) \mapsto (\epsilon_1x_1,-v/(u^2-v^2),u/(u^2-v^2)), \\
\tau &:& (x_1,u,v) \mapsto (1/x_1,v,u), \\
\lambda &:& (x_1,u,v) \mapsto (a/x_1,\epsilon_2\epsilon_3u,\epsilon_2\epsilon_3v).
\end{eqnarray*}

Define $y_2=u+v$, $y_3=u-v$. Then we get
\begin{eqnarray*}
\sigma &:& (y_2,y_3) \mapsto (1/y_2,-1/y_3), \\
\tau &:& (y_2,y_3) \mapsto (y_2,-y_3), \\
\lambda &:& (y_2,y_3) \mapsto (\epsilon_2\epsilon_3y_2,\epsilon_2\epsilon_3y_3).
\end{eqnarray*}

Thus we may apply Theorem \ref{t3.6} and conclude that $K(x_1,y_2,y_3)^{\langle\sigma,\tau,\lambda\rangle}$
is rational over $K$.
\end{idef}

\bigskip
\begin{idef}{Case 2.} $G=W_2(195)$.

The action of $G=\langle\sigma,\tau,\lambda\rangle$ is given by
\begin{eqnarray*}
\sigma &:& (x_1,x_2,x_3)\mapsto (\epsilon x_1,x_3,ax_1/x_2), \\
\tau &:& (x_1,x_2,x_3)\mapsto (1/x_1,b/x_3,b/x_2), \\
\lambda &:& (x_1,x_2,x_3) \mapsto (c_1/x_1,c_2/x_2,c_3/x_3)
\end{eqnarray*}
where $a,b,c_1,c_2,c_3\in K\backslash \{0\}$, $\epsilon\in\{1,-1\}$
(see Case 5. $G=W_{11}(187)$ in the proof of Theorem \ref{t3.7}).

Define $y_1=x_1$, $y_2=x_2/\sqrt{b}$, $y_3=\sqrt{b}/x_3$ and use the relations $\tau\sigma\tau^{-1}=\sigma^{-1}$,
$\sigma\lambda=\lambda\sigma$, $\tau\lambda=\lambda\tau$.
We find that
\begin{eqnarray*}
\sigma &:& (y_1,y_2,y_3)\mapsto (\epsilon y_1,1/y_3,ay_2/y_1), \\
\tau &:& (y_1,y_2,y_3)\mapsto (1/y_1,y_3,y_2), \\
\lambda &:& (y_1,y_2,y_3)\mapsto (\epsilon_1/y_1,\epsilon_2/y_2,\epsilon_2/y_3)
\end{eqnarray*}
where $\epsilon,\epsilon_1,\epsilon_2\in \{1,-1\}$ and $a^2=\epsilon_1$.

We will ``copy" the proof of Case 5. $G=W_{11}(187)$ in Theorem
\ref{t3.7} without changing anything. Note that
$\lambda:(u,v)\mapsto(\epsilon\epsilon_2au/y_1,\epsilon\epsilon_2avy_1)$.

Following the proof of $G=W_{11}(187)$, we define $z_2=u+y_1v$,
$z_3=u-y_1v$. Since $\lambda(y_1)=\epsilon_1/y_1$, we find that
$\lambda:(z_2,z_3)\mapsto(\epsilon\epsilon_2az_2/y_1,\epsilon\epsilon_2az_3/y_1)$
if $\epsilon_1=1$; while
$\lambda:(z_2,z_3)\mapsto(\epsilon\epsilon_2az_3/y_1,\epsilon\epsilon_2az_2/y_1)$
if $\epsilon_1=-1$.

Thus we may discuss various cases when $\epsilon,\epsilon_1=\pm 1$ and prove that
$K(y_1,z_2,z_3)^{\langle\sigma,\tau,\lambda\rangle}$ is rational over $K$ by Theorem \ref{t3.6}.
\end{idef}
\end{proof}

\begin{proof}[\indent Proof of Theorem \ref{t1.3}.]
By all the theorems proved in this section, i.e. Theorems
\ref{t3.1} -- \ref{t3.8}, the proof of Theorem \ref{t1.3} is
completed.
\end{proof}

\setcounter{equation}{0}
\section{Another proof of Theorem \ref{t1.3}}

In this section we will give an alternative, short proof of
Theorem \ref{t1.3} in the case where the ground field $K$ is
algebraically closed and has characteristic $0$. The proof is
geometric and does not require computations. Moreover it does not
use the classification of finite subgroups in $GL_3(\bm{Z})$
\cite{Ta}.

Throughout this section we assume that $G$ is a finite $2$-group
and $K$ is an algebraically closed field of characteristic $0$.

The following Lemma \ref{l4.1}, Lemma \ref{l4.3}, Corollary
\ref{c4.2} and Corollary \ref{c4.4} are meant to replace a
classification of finite subgroups of $GL_2(\bm{Z})$ and
$GL_3(\bm{Z})$. (In fact, finite subgroups of $GL_2(\bm{Z})$ are
listed in \cite{Ha1,Ha2} and those of $GL_3(\bm{Z})$ are listed in
\cite{Ta}.)

Lemma \ref{l4.1} is an easy exercise in the representation theory
of finite groups.

\begin{lemma} \label{l4.1}
Let $G$ be a finite $2$-group. Then any representation $G\to
GL_3(\bm{Q})$ is reducible.
\end{lemma}

\begin{coro} \label{c4.2}
Let $G$ be a finite $2$-group. Then any
integral representation $G\to GL_3(\bm{Z})$ has a subrepresentation
of $\bm{Z}$-rank $1$.
\end{coro}

Now let $G$ be a finite $2$-group acting monomially on $K(x,y,z)$.
Let $\rho:G\to GL_3(\bm{Z})$ be the integral representation
introduced in Definition \ref{d2.7}. By Corollary \ref{c4.2} there
is a $1$-dimensional subrepresentation $\rho_1$. Denote the
quotient representation by $\rho_2 : G\to GL_2(\bm{Z})$. In a
suitable basis in $\bm{Z}^3$ any element of $\rho(G)$ can be
written as
\[
\rho(g)=
\begin{pmatrix}
a & b_1&b_2 \\
0& c_{11}& c_{12} \\
0&c_{21}& c_{22}
\end{pmatrix}
\]
and then
\[
\rho_2(g)=\begin{pmatrix} c_{11}& c_{12} \\ c_{21}& c_{22} \end{pmatrix}.
\]

We will denote by $M$ the lattice $\bm{Z}^3$. Thus the group $G$
acts on $M$. By Lemma \ref{l2.8} we may assume that the
representation $\rho$ is faithful on $M$.

\begin{lemma} \label{l4.3}
Any $2$-subgroup in $GL_2(\bm{Z})$ is
conjugated to a subgroup of the group generated by two matrices
$A$ and $B$ $($and isomorphic to $D_4)$
\[
A= \begin{pmatrix} 0&-1 \\ 1& 0 \end{pmatrix},
\qquad B= \begin{pmatrix} 1&0 \\ 0& -1 \end{pmatrix}.
\]
\end{lemma}

\begin{proof}[\indent Outline of the proof]
For any matrix $C$ in the group we have $\fn{tr} C=0$ or $\pm 1$
and $\det C=\pm 1$. Since $C$ is an element of order $2^k$, its
minimum polynomial divides $t^{2^k}-1$. This immediately implies
$\fn{tr} C=0$ and $C^4=E$. Hence $C$ is conjugate to either $A$ or
$B$, or $\pm E$. On the other hand, the order of the group is at
most $8$ by the famous Minkowski theorem (see, for example,
\cite[page 169]{Ta}).
\end{proof}

Let $N\subset G$ be the maximal subgroup that acts trivially on
$L$. It is easy to see that so is the restriction $\rho_2|_N$.

\begin{coro} \label{c4.4}
In the above notation, there are only the following possibilities:
\begin{equation}
\label{eq4.1} N\simeq D_4,\quad C_4, \quad C_2, \quad
\text{or}\quad N=\{1\}.
\end{equation}
\end{coro}

By Lemma \ref{l4.3}, for any $g\in G$, the action of $g$ on
$K(x,y,z)$ can be written as follows:
\begin{eqnarray*}
g_1&:& x \longmapsto \alpha_1 x^{a_1},~ y \longmapsto \Psi_1 z, ~z \longmapsto \Psi_2 y^{-1}, \\
g_2&:& x \longmapsto \alpha_2 x^{a_2},~ y \longmapsto \Psi_3 y, ~z \longmapsto \Psi_4 z^{-1},
\end{eqnarray*}
where $a_i\in \{ 1, -1 \}$, $\alpha_i\in K^{\times}$, and
$\Psi_i\in K(x)^{\times}$ are some monomials, i.e. elements of the
form $cx^n$ where $c \in K^{\times}$ and $n$ is an integer.

\begin{remark} \label{r4.6} \rm
Define $L:=K(x)$. We may regard $K(x,y,z)$ as the function field
of the projective surface $\bm{P}^1_L\times \bm{P}^1_L$ over $L$
and the induced action of $G$ on $\bm{P}^1_L\times \bm{P}^1_L$ is
semi-linear. Here $y$ and $z$ are non-homogeneous coordinates on
the first and second factors respectively.
\end{remark}

We use some facts about del Pezzo surfaces. Good references are
\cite{Manin-1966,Manin-Cubic-forms-e-II,0621.14029}. Recall that
over an algebraically closed field there are two types of del
Pezzo surfaces of degree $8$: the Hirzebruch surface $\bm{F}_1$
and $\bm{P}^1\times \bm{P}^1$. A form $S$ of $\bm{F}_1$ over any
field $L$ is not minimal, that is, it contains a unique
$(-1)$-curve. Contracting it we get a Brauer-Severi variety $S'$
having an $L$-point. Hence, $S'\simeq \bm{P}^2_L$. A form of
$\bm{P}^1\times \bm{P}^1$ over a field $L$ is called a
\textit{minimal del Pezzo surface of degree $8$}. Taking into
account the remark after Corollary \ref{c4.4}, one find that
Theorem \ref{t1.3} is a consequence of Proposition \ref{p4.7}
below.

\begin{proposition} \label{p4.7}
Let $L$ be a field of characteristic $0$ and $G$ be a finite
$2$-group. Let $S$ a minimal del Pezzo surface of degree $8$
defined over $L$. Suppose $G$ acts on both $L$ and $S$ which is
regarded as an $L^G$-scheme, and let $N$ be the kernel of the
homomorphism $G\to \operatorname{Aut} L$. Assume that $N$ is given
as in \eqref{eq4.1}. If $L^G$ is a $C_1$-field, then $S/G$ is
rational over $L^G$.
\end{proposition}

For the definition $C_i$-fields, see \cite[page 3; DKT, page
1256]{Gr}. An algebraic function field of one variable over an
algebraically closed field is a $C_1$-field by Tsen's Theorem
\cite[page 22; DKT, Theorem 5.8, page 1259]{Gr}.

\begin{proof}
Our proof is by induction on the order of $G$.

If $G=\{1\}$, then $S$ has an $L$-point (because $L$ is a $C_1$-
field) and so $S$ is a rational quadric in $\bm{P}^3$.

From now on we assume that $n:=|G|>1$ and the assertion holds for
all groups $G_1$ with $|G_1|<n$.

Let $\bar L$ be the algebraic closure. Set $\bar S:=S\otimes \bar
L$. The group $N$  acts linearly on $S$ over $L$ and on $\bar S$
over $\bar L$. In particular, $N$ acts on $\fn{Pic} (\bar S)\simeq
\bm{Z}\oplus \bm{Z}$.

First we assume that $|N|>2$. By Formula \eqref{eq4.1} there is an
element $g\in N$ of order $4$ and $\tau:=g^2$ is contained in
$Z(G)\cap N$, where $Z(G)$ is the center of $G$. Indeed, in this
case, $\tau$ is the only element of order $2$ in $Z(N)\supset
Z(G)\cap N$. Denote  $N_0:=\{1,\tau\}$ and $F:=S/N_0$. By our
assumptions, $F$ is defined over $L$ and $G/N_0$ acts naturally on
$F$ by semi-linear automorphisms. Thus $S/G\simeq F/(G/N_0)$. In a
suitable (non-homogeneous) coordinate system in $\bar F\simeq
\bm{P}^1_{\bar L}\times \bm{P}^1_{\bar L}$
 over $\bar L$ there are only the following
possibilities for the action:
\begin{enumerate} \itemsep=-2pt
\item[(1)] $\tau : x \longmapsto y,\quad y\longmapsto x$;
\item[(2)] $\tau : x \longmapsto x,\quad y \longmapsto -y$;
\item[(3)] $\tau : x \longmapsto -x,\quad y \longmapsto -y$.
\end{enumerate}
The element $\tau=g^2$  acts trivially on $\fn{Pic} (\bar S)\simeq
\bm{Z}\oplus \bm{Z}$. Thus the case (1) is impossible. In the case
(2) we have $\bar F\simeq \bm{P}^1_{\bar L}\times \bm{P}^1_{\bar
L}$, so $F$ a minimal del Pezzo surface of degree $8$. We reduce
the problem to an action of a smaller group $G_1=G/N$ on $F$.

Consider the case (3). Let $\pi : S\to F$ be the quotient
morphism. Over the field $\bar L$, the group $N$ acts freely in
codimension one on $S$ and has exactly 4 fixed points. These
points give us 4 singular points on $F$ of type $A_1$ and $F$  has
no other singularities. In particular, this implies that the
divisor $K_F$ is Cartier. Since the map $\pi: S\to F$ is \'etale
in codimension one, $K_S=\pi^* K_F$. Hence $-K_{F}$ is  ample and
$K_{F}^2=K_S^2/2=4$. This means that $F$ is a del Pezzo surface of
degree 4 having 4 singular points of type $A_1$. According to the
classification \cite{Coray1988}, there are exactly 4 lines
$l_1,\dots,l_4$ on $F$. Let $\mu :\tilde F \to F$ be the minimal
resolution and let $\tilde l_i$ be the proper transform of $l_i$.
Then $\tilde F$ is a \textit{weak del Pezzo} surface in the sense
that its anti-canonical divisor $-K_{\tilde F}$ is nef and big
(such surfaces are also called generalized del Pezzo surfaces
\cite{Coray1988}). Moreover, $K_{\tilde F}^2=K_F^2=4$. Again by
\cite[Prop. 6.1]{Coray1988} $\tilde l_1,\dots,\tilde l_4$ are
disjointed $(-1)$-curves (and $\tilde F$ contains no other
$(-1)$-curves). Consider their contraction $\varphi : \tilde F\to
F'$. Here $K_{F'}^2=K_{\tilde F}^2+4=8$. Further,
$-K_{F'}=\varphi_*(-K_{\tilde F})$ is nef and big, i.e., $F'$ is a
weak del Pezzo surface of degree 8. Moreover, every $(-2)$-curve
$C$ on $\tilde F$ meets two $\varphi$-exceptional curves  $\tilde
l_i$. Hence the surface $F'$ contains no curves with negative
self-intersection numbers. Therefore, $F'$ is a minimal  del Pezzo
surface of degree 8. Finally, both the birational maps $\mu :
\tilde F \rightarrow F'$ and $\varphi : \tilde F \rightarrow F$
are canonically defined, so they are $G/N$-equivariant. Again we
reduce the problem to an action of $G_1=G/N$ on $F'$ with $|G_1| <
|G|$.

Now assume that $|N|=2$, so $N=\{1,\tau\}$. As above, $S/G\simeq
F/(G/N_0)$ and in a suitable coordinate system in $\bar F\simeq
\bm{P}^1_{\bar L}\times \bm{P}^1_{\bar L}$ over $\bar L$ there are
only the possibilities (1), (2) and (3) listed above. In cases (2)
and (3) we may argue as above. It remains to consider the case
(1). Then the quotient $\bar F=\bar S/ N$ is isomorphic to
$\bm{P}^2_{\bar L}$. Thus $F$ is a del Pezzo surface  of degree
$9$ (a Brauer-Severi variety). In this situation, $G_1=G/N$ acts
on $L$ effectively. Then by Lemma \ref{l4.8} below the quotient $S
/ G\simeq \bm{P}^2_L/(G/N)$ is $L$-rational.

Finally, if $N=\{1\}$, i.e., the action of $G$ on $L$ is
effective, we may also apply Lemma \ref{l4.8} to finish the proof.
\end{proof}

\begin{lemma} \label{l4.8}
Let $\Bbbk$ be a field of characteristic $0$ and let $S$ be a del
Pezzo surface of degree $d\ge 5$ over $\Bbbk$. Suppose a finite
group $G$ acts effectively  on both $\Bbbk$ and $S$ which is
regarded as a $\Bbbk^G$-scheme. If $\Bbbk^G$ is a $C_1$-field,
then $S/G$ is rational over $\Bbbk^G$.
\end{lemma}

\begin{proof}
We claim that the natural $G$-equivariant morphism $S \rightarrow
S/G \times_{\Bbbk^G} \Bbbk$ is an isomorphism. Since it is a local
property, it suffices to show that, if $G$ acts on an affine
domain $R$ over a field $\Bbbk$ so that $\Bbbk$ is left invariant
by the action of $G$ and $G$ acts faithfully both on $\Bbbk$ and
$R$, then $R^G\otimes_{\Bbbk^G}\Bbbk \simeq R$. This follows from
the Galois descent lemma \cite[Lemma 18.1, page 279]{KMRT}. Note
that the Galois descent lemma is essentially another form of
Theorem \ref{t2.1}.

Now consider the following diagram
\[
\begin{array}{c@{}c@{}c}
S & \longrightarrow & S/G \\[2pt]
\big\downarrow & & \big\downarrow \\[2pt]
\operatorname{Spec} \Bbbk & \to & \operatorname{Spec} \Bbbk^G
\end{array}
\]
Since $S/G \times_{\Bbbk^G} \Bbbk\simeq S$, it follows that $S/G$
is a minimal del Pezzo surface of degree $d\ge 5$. By
\cite[Chapter 4, \S7]{Manin-Cubic-forms-e-II} it is sufficient to
show that $S/G$ has a $\Bbbk^G$-point. This holds for an arbitrary
$\Bbbk$ if $d=5$ \cite{Swinnerton-Dyer1972} and for $C_1$-fields
if $d\ge 6$ \cite{Manin-1966}. Thus $S/G$ is $\Bbbk^G$-rational.
\end{proof}

\section{Applications}

\begin{proof}[\indent Proof of Theorem \ref{t1.4}]

By \ref{t2.2}, we may assume that $V$ contains no $1$-dimensional
subrepresentation. By \cite[pages 53 and 65]{Se}, every
irreducible representation of $G$ is an induced representation
from some 1-dimensional representation of $H$ for some subgroup
$H$ of $G$; thus it is a monomial representation. It follows that
$V$ is isomorphic to (i) an irreducible 4-dimensional
representation, (ii) the direct sum of two irreducible
2-dimensional representations, or (iii) the direct sum of three
2-dimensional irreducible representations.

\begin{idef}{Case 1.} $V$ is an irreducible 4-dimensional representation.

Since $V$ is induced from a 1-dimensional representation,
for any $\sigma\in G$, $\sigma(e_i)=a_i(\sigma)\cdot e_j$ for some basis $e_1$, $e_2$, $e_3$, $e_4$ of $V$ where
$a_i(\sigma)\in K\backslash\{0\}$.
If $x_1$, $x_2$, $x_3$, $x_4$ denotes the dual basis of $e_1$, $e_2$, $e_3$, $e_4$,
then $K(V)=K(x_1,x_2,x_3,x_4)$ and $\sigma\cdot x_i=b_i(\sigma)\cdot x_j$ for any $\sigma\in G$ where
$b_i(\sigma)\in K\backslash \{0\}$.

Define $y_i=x_i/x_4$ for $1\le i\le 3$. Then
$\sigma(x_4)=t_\sigma\cdot x_4$ where $t_\sigma \in
K(y_1,y_2,y_3)$ for any $\sigma\in G$. Apply Theorem \ref{t2.2}.
We find that $K(x_1,x_2,x_3,x_4)^G=K(y_1,y_2,y_3)^G(f)$ for some
$f$ satisfying $\sigma(f)=f$ for any $\sigma\in G$. Note that $G$
acts on $K(y_1,y_2,y_3)$ by monomial $K$-automorphisms. Apply
Theorem \ref{t1.3}. We find that $K(y_1,y_2,y_3)^G$ is rational
over $K$. Thus the quotient $P(V)/G$ is also rational.
\end{idef}

\bigskip
\begin{idef}{Case 2.} $V$ is a direct sum of irreducible 2-dimensional representations and
certain 1-dimensional representations.

Apply Theorem \ref{t2.2} and Theorem \ref{t1.1}. The details are
omitted.
\end{idef}

\bigskip
\begin{idef}{Case 3.} $V\simeq V_1\oplus V_2\oplus V_3$ where each
$V_i$ is the representation space of an irreducible
$2$-dimensional representation of $G$. It follows that
$K(V)=K(x_1,x_2,y_1,y_2,z_1,z_2)$ and $G$ acts on elements of the
set $\{ x_1, x_2 \}$ (resp. $\{ y_1, y_2 \}$, $\{ z_1, z_2 \}$) by
permutations up to non-zero elements in $K$. Define $x=x_1/x_2,
y=y_1/y_2, z=z_1/z_2$. Apply Theorem 2.2 repeatedly. We reduce the
question to the rationality of $K(x,y,z)^G$. Apply Theorem 1.3.
\end{idef}
\end{proof}

For the proof of Theorem \ref{t1.5}, we need two standard results
in group theory, whose proof will be omitted.

\begin{lemma} \label{l5.2}
Let $p$ be a prime number, $G$ be a finite $p$-group.
Denote by $Z(G)$ the center of $G$.
Let $H$ be a normal subgroup of $G$.
If $H\supsetneq \{1\}$, then $H\cap Z(G) \supsetneq \{1\}$.
\end{lemma}

\begin{lemma} \label{l5.3}
Let $G$ be a finite $2$-group of exponent $e$, and $K$ be a field
containing a primitive $e$-th root of unity. For each non-negative
integer $i$, let $s_i$ be the number of non-equivalent
$2^i$-dimensional irreducible representations of $G$ over $K$.
Then $\sum_{i}s_i4^i = |G|$.
\end{lemma}

\begin{proof}[\indent Proof of Theorem \ref{t1.5} when $K$ is algebraically closed.]

Let $G$ be a group of order 32. If $G$ is abelian, then $K(G)$ is
rational by Fischer's Theorem \cite[Theorem 6.1]{Sw}. Thus we will
assume that $G$ is non-abelian from now on. Let $Z(G)$ be the
center of $G$, and $[G,G]$ be the commutator subgroup of $G$.

If $\fn{char}K=2$, then $K(G)$ is rational by Kuniyoshi's Theorem \cite{Ku}.
Thus we will assume that $\fn{char}K\ne 2$ and $K$ is algebraically closed from now on.

\bigskip
\begin{idef}{Case 1.} $|Z(G)|\ge 8$.

Since $|G/Z(G)|\le 4$, $G/Z(G)$ is isomorphic to $C_2\times C_2$.
Thus there is an abelian subgroup $H$ so that $Z(G)\subset
H\subset G$ and $[G:H]=2$. Apply Theorem \ref{t2.6}. We find that
$K(G)$ is rational.
\end{idef}

\bigskip
\begin{idef}{Case 2.} $Z(G)\simeq C_2$ or $C_4$.

Let $\tau\in Z(G)$ be the unique element of order 2 in $Z(G)$.

We claim that there is an irreducible representation $\rho:G\to
GL(V)$ over $K$, which is faithful.

Assume the validity of this claim. Note that $\dim_K V\le 4$ by
Lemma \ref{l5.3}. Hence $K(V)^G$ is rational over $K$ by Theorem
\ref{t1.4}. Since $K(G)$ is rational over $K(V)^G$ by Theorem
\ref{t2.1}, we find that $K(G)$ is rational.

Now we will prove this claim. Suppose not, i.e. every irreducible
representation $\rho:G\to GL(W)$ is not faithful. Then
$\fn{Ker}(\rho)\cap Z(G) \supsetneq \{1\}$ by Lemma \ref{l5.2}.
Thus $\tau\in\fn{Ker}(\rho)\cap Z(G)$. We conclude that
$\rho(\tau)$ is trivial for all irreducible representations of
$G$. But the regular representation of $G$ is faithful and is the
direct sum of these irreducible representations of $G$. This leads
to a contradiction.
\end{idef}

\bigskip
\begin{idef}{Case 3.} $Z(G)\simeq C_2\times C_2$.

We will show that $G$ has a faithful 4-dimensional representation,
which is either irreducible or a direct sum of two 2-dimensional irreducible representations.
In any case, the rationality of $K(G)$ can be proved as in Case 2.

\bigskip
Step 1. If $G$ has a faithful 4-dimensional representation,
we are done. From now on, we will assume that, if $G$ has a 4-dimensional irreducible representation,
then it is not faithful.

We will prove that $G$ has no 4-dimensional irreducible representation at all.

Suppose that $\rho:G\to GL_4(K)$ is irreducible. Then
$\fn{Ker}(\rho)\supsetneq \{1\}$. Thus $\rho$ gives rise to a
faithful irreducible 4-dimensional representation
$G/\fn{Ker}(\rho)\to GL_4(K)$. Note that $|G/\fn{Ker}(\rho)|\le
16$ which is impossible by Lemma \ref{l5.3}.

\bigskip
Step 2. We will show that $G$ has at least four non-equivalent
2-dimensional irreducible representation.

Let $G$ has $p$ 1-dimensional irreducible representations and has
$s$ 2-dimensional irreducible representations. Then $p+4s=|G|=32$
and $p=|G/[G,G]|$ which is a divisor of 16. The only solutions of
$(p,s)$ are $(p,s)=(4,7)$, $(8,6)$, $(16,4)$. Thus $s\ge 4$.

\bigskip
Step 3. From Step 2, let $\rho_i:G\to GL_2(K)$ be all the 2-dimensional irreducible representations
where $1\le i\le s$ with $s\ge 4$.

We claim that, there is some $\rho_i$ so that $\fn{Ker}(\rho_i) \not\supset Z(G)$.

Otherwise, assume that $\fn{Ker}(\rho_i)\supset Z(G)$ for all
$1\le i\le s$. Note that the kernel of any 1-dimensional
representation contains $[G,G]$. Hence $[G,G]\cap (\bigcap_{1\le
i\le s}\fn{Ker}(\rho_i)) \supset [G,G] \cap Z(G) \supsetneq \{1\}$
by Lemma \ref{l5.2}. It follows that there is some element $g\in
G$, $g\ne 1$, but $g$ belongs to the kernel of any irreducible
representation of $G$. This is impossible by the same arguments as
in the last paragraph of Case 2.

\bigskip
Step 4. By Step 3, we may assume that $\fn{Ker}(\rho_1)\not\supset
Z(G)$. Write $\fn{Ker}(\rho_1)\cap Z(G)=\langle\lambda\rangle$
where $\lambda$ is an element of order two in $Z(G)\simeq
C_2\times C_2$.

Define $H=\bigcap_{1\le i\le s} \fn{Ker}(\rho_i)$. We will prove that $H=\{1\}$.

Suppose that $H\supsetneq \{1\}$. But Lemma \ref{l5.2}, $H\cap
Z(G) \supsetneq \{1\}$. But $H\cap Z(G)\subset
\fn{Ker}(\rho_1)\cap Z(G)=\langle\lambda\rangle$. Hence $H\cap
Z(G)=\langle\lambda\rangle$. Thus the group
$G/\langle\lambda\rangle$ has at least four 2-dimensional
irreducible representations. If the number of 1-dimensional
irreducible representations of $G/\langle\lambda\rangle$ is $q$,
then $q+4\cdot 4 \le |G/\langle\lambda\rangle|$ by Lemma
\ref{l5.3}. Since $|G/\langle\lambda\rangle|=16$, we get a
contradiction.

\bigskip
Step 5. Since $H=\{1\}$, we find that $\lambda\notin H$.
We may assume that $\lambda\notin \fn{Ker}(\rho_2)$.

Define $H_0=\fn{Ker}(\rho_1)\cap \fn{Ker}(\rho_2)$. If
$H_0\supsetneq \{1\}$, then $H_0\cap Z(G)\supsetneq \{1\}$ by
Lemma \ref{l5.2}. Hence $H_0\cap Z(G)\subset \fn{Ker}(\rho_1)\cap
Z(G)=\langle\lambda\rangle$. Thus $H_0\cap
Z(G)=\langle\lambda\rangle$, which is impossible because
$\lambda\notin \fn{Ker}(\rho_2)$.

\bigskip
Step 6. Since $\{1\}=H_0=\fn{Ker}(\rho_1)\cap\fn{Ker}(\rho_2)$,
the direct sum of $\rho_1$ and $\rho_2$ is a faithful representation of $G$. Done.
\end{idef}
\end{proof}

\begin{remark} \em
As we note in Section 1, the above proof is valid only under the
assumption that $\sqrt{a}\in K$ for any $a\in K$ instead of the
weaker assumption that $\zeta_e\in K$ (which is given in Theorem
\ref{t1.5}). Consider the following proposition : \textit{Let G be
a non-abelian group of order $32$ and let $e$ be the exponent of
$G$. If $K$ is a field containing $\zeta_e$ and $G \rightarrow
GL(V)$ is an irreducible representation over $K$ where $dim V=4$,
then $K(V)^G$ is rational over $K$}. Note that the above
proposition is valid if we assume the classification of groups of
order $32$, because there are exactly seven such groups (they are
groups $(42;49), (43;50), (44;43), (45;44), (46;6), (47;7),
(48;8)$ in the notation of \cite[page 3025]{CHKP}) and all the
$4$-dimensional irreducible representations can be described
explicitly. On the other hand, at present we don't have a proof of
the above proposition if we try to avoid using the classification
of groups of order $32$. With this proposition in hand, the above
proof of Theorem \ref{t1.5} may be adapted to the situation
assuming only that $\zeta_e\in K$.
\end{remark}

\newpage
\renewcommand{\refname}{\centering{References}}


\begin{thebibliography}{BEO}

\bibitem[AHK]{AHK}
H. Ahmad, M. Hajja and M. Kang,
\textit{Rationality of some projective linear actions},
J. Algebra 228 (2000), 643--658.

\bibitem[CHK]{CHK}
H. Chu, S.-J. Hu and M. Kang, \textit{Noether's problem for
dihedral 2-groups}, Comment. Math. Helv. 79 (2004), 147--159.

\bibitem[CHKK]{CHKK}
H. Chu, S.-J. Hu, M. Kang and B. E. Kunyavskii, \textit{Noether's
problem and the unramified Brauer group for groups of order 64},
preprint.

\bibitem[CHKP]{CHKP}
H. Chu, S.-J. Hu, M. Kang and Y. G. Prokhorov, \textit{Noether's
problem for groups of order $32$}, J. Algebra 320 (2008),
3022--3035.

\bibitem[CT]{Coray1988}
D. F. Coray and M. A. Tsfasman, \textit{Arithmetic on singular del
Pezzo surfaces}, Proc. Lond. Math. Soc. 57 (1988) 25--87.

\bibitem[DKT]{DKT}
S. Ding, M. Kang and E. Tan \textit{Chiungtze C. Tsen (1898 -
1940) and Tsen's theorems}, Rocky Mount. J. Math. 29 (1999),
1237--1269.

\bibitem[Gr]{Gr}
M. J. Greenberg, \textit{Lectures on forms in many variables}, W.
A. Benjamin, Inc., New York, 1969.

\bibitem[Ha1]{Ha1}
M. Hajja, \textit{A note on monomial automorphisms}, J. Algebra 85
(1983), 243--250.

\bibitem[Ha2]{Ha2}
M. Hajja,
\textit{Rationality of finite groups of monomial automorphisms of $K(x,y)$},
J. Algebra 109 (1987), 46--51.

\bibitem[HK1]{HK1}
M. Hajja and M. Kang,
\textit{Finite group actions on rational function fields},
J. Algebra 149 (1992), 139--154.

\bibitem[HK2]{HK2}
M. Hajja and M. Kang,
\textit{Three-dimensional purely monomial group actions},
J. Algebra 170 (1994), 805--860.

\bibitem[HK3]{HK3}
M. Hajja and M. Kang,
\textit{Some actions of symmetric groups},
J. Algebra 177 (1995), 511--535.

\bibitem[HKO]{HKO}
M. Hajja, M. Kang and J. Ohm, \textit{Function fields of conics as
invariant subfields}, J. Algebra 163 (1994), 383--403.


\bibitem[HR]{HR}
A. Hoshi and Y. Rikuna, \textit{Rationality problem of
three-dimensional purely monomial group actions: the last case},
Math. Comp. 77 (2008), 1823--1829.

\bibitem[Is]{Is}
V. A. Iskovskikh,
\textit{Rational surfaces with a pencil of rational curves},
Math. USSR Sb. 3 (1967), 563--587.

\bibitem[Ka1]{Ka1}
M. Kang, \textit{Rationality problem of $GL_4$ group actions},
Advances in Math. 181 (2004), 321--352.

\bibitem[Ka2]{Ka2}
M. Kang, \textit{Some group actions on $K(x_1,x_2,x_3)$}, Israel
J. Math. 146 (2005), 77--93.

\bibitem[Ka3]{Ka3}
M. Kang, \textit{Some rationality problems revisited}, in
``Proceedings of the 4th International Congress of Chinese
Mathematicians, Hangzhou, 2007", edited by L. Ji, Kefeng Liu, Lo
Yang and Shing-Tung Yau, Higher Education Press (Beijing) and
International Press (Somerville).

\bibitem[Ka4]{Ka4}
M. Kang, \textit{Rationality problem for some meta-abelian
groups}, to appear in ``J. Algebra".

\bibitem[KMRT]{KMRT}
M.-A. Knus, A. Merkurjev, M. Rost and J.-P. Tignol, \textit{The
book of involutions}, American Mathematical Society, Providence,
1998.


\bibitem[Ku]{Ku}
H. Kuniyoshi,
\textit{On a problem of Chevalley}, Nagoya Math. J. 8 (1955), 65--87.

\bibitem[Ma1]{Manin-1966}
Yu. I. Manin. \textit{Rational surfaces over perfect fields},
Inst. Hautes \'Etudes Sci. Publ. Math. 30 (1966) 55--113.

\bibitem[Ma2]{Manin-Cubic-forms-e-II}
Yu. I. Manin, \textit{Cubic forms  {\rm : } Algebra, geometry and
arithmetic}, second edition, English translation by M. Hazewinkel,
North-Holland Publ. Co., Amsterdam,
 1986.

\bibitem[MT]{0621.14029}
Yu.I. Manin and M.A. Tsfasman, \textit{Rational varieties {\rm
:}Algebra, geometry and arithmetic}, Russ. Math. Surv. 41 (1986)
51--116.

\bibitem[Pr]{Pr}
Y. G. Prokhorov, \textit{Fields of invariants for finite linear
groups}, in ``Rationality problem" edited by F. Bogomolov and Y.
Tschinkel, Progress in Math., Birkhauser, Boston, to appear.

\bibitem[Sa]{Sa}
D. J. Saltman, \textit{A nonrational field, answering a question
of Hajja}, in ``Algebra and Number Theory" edited by M. Boulagouaz
and J.-P. Tignol, Marcel Dekker, New York, 2000.

\bibitem[SD]{Swinnerton-Dyer1972}
H.P.F. Swinnerton-Dyer, \textit{Rational points on del Pezzo
surfaces of degree 5}, in ``Algebraic Geom., Oslo 1970, Proc. 5th
Nordic Summer-School Math.", Wolters-Noordhoff, Groningen, 1972.


\bibitem[Se]{Se}
J.-P. Serre,
\textit{Linear representations of finite groups},
Springer GTM vol. 42, Springer-Verlag, 1978.

\bibitem[Sw]{Sw}
R. G. Swan,
\textit{Noether's problem in Galois theory},
in ``Emmy Noether in Bryn Mawr", edited by B. Srinivasan and J. Sally,
Springer-Verlag, Berlin, 1983.

\bibitem[Ta]{Ta}
K. Tahara,
\textit{On the finite subgroups of $GL(3,\bm{Z})$},
Nagoya Math. J. 41 (1971), 169--209.





\end{thebibliography}
\end{document}